\documentstyle[12pt]{amsart}

\newtheorem{theorem}{THEOREM}[section]
\newtheorem{definition}[theorem]{Definition}
\newtheorem{question}{Question}

\newtheorem{corollary}[theorem]{Corollary}
\newtheorem{lemma}[theorem]{Lemma}
\newtheorem{example}{EXAMPLE}

\newtheorem{proposition}[theorem]{Proposition}
\def\hbar{\overline{h}}
\def \bomeg{\partial \Omega}
\def\Omegabar{\overline{\Omega}}

\def\endpf{\qed}

\def\epf{\hskip.2in\vrule width.4pt height6.65pt
depth.15pt\vrule
width2.5pt height6.65pt depth-6.25pt\hskip-2.5pt\vrule
width2.5pt
height.25pt depth.15pt\vrule width.4pt
height6.65pt depth.15pt\ }

\def\proof{\noindent {\bf Proof. }}

\def\qed{\hskip.2in\vrule width.4pt
height6.65pt depth.15pt\vrule
width2.5pt height6.65pt depth-6.25pt\hskip-2.5pt\vrule
width2.5pt
height.25pt depth.15pt\vrule width.4pt
height6.65pt depth.15pt\ }

\def \Omegabar{\overline \Omega}
\def \bomega{\partial \Omega}

\def \d{\partial}

\def \X{\cal{ X} }

\def \P{\cal P }
\def\zbar{\overline{z}}

\def \hbar{\overline{h}}

\def\H {\cal H}

\def\supp{\hbox{supp}}

\def\p{\bf p}
\font\teneufm=eufm10
\font\seveneufm=eufm7
\font\fiveeufm=eufm5
\newfam\eufmfam
\textfont\eufmfam=\teneufm
\scriptfont\eufmfam=\seveneufm
\scriptscriptfont\eufmfam=\fiveeufm

\newfam\msbfam
\font\tenmsb=msbm10    \textfont\msbfam=\tenmsb
\font\sevenmsb=msbm7   \scriptfont\msbfam=\sevenmsb
\font\fivemsb=msbm5    \scriptscriptfont\msbfam=\fivemsb
\def\Bbb{\fam\msbfam \tenmsb}

\def\RR{{\Bbb R}}
\def\CC{{\Bbb C}}

\def\HH{{\Bbb H}}

\def\T{\cal T}
\newfam\bigfam
\font\tenbig=msbm10 scaled \magstep2   \textfont\bigfam=\tenbig
\font\sevenbig=msbm7 scaled \magstep2  \scriptfont\bigfam=\sevenbig
\font\fivebig=msbm5 scaled \magstep2   \scriptscriptfont\bigfam=\fivebig

\begin{document}

\title[Hardy Classes, Integral Operators, and Duality]
{Hardy Classes, Integral Operators, and Duality on
 Spaces of Homogeneous Type}

\author{Steven G. Krantz}
\author{Song-Ying Li}
\address{\hskip-\parindent
Steven G. Krantz\\
Department of Mathematics\\
Washington University\\
Such that. Louis, MO 63130}
\email{sk@@math.wustl.edu}
\address{\hskip-\parindent
Song-Ying Li\\
Department of Mathematics\\
Washington University\\
Such that. Louis, MO 63130}
\email{songying@@math.wustl.edu}

\thanks{Krantz's research is supported in part
by NSF Grant DMS-9022140 during his stay at MSRI.}

\begin{abstract}
 The authors study Hardy spaces,
of arbitrary order, on a space
of homogeneous type.  This extends earlier work that treated only
$H^p$ for $p$ near 1.  Applications are given to the boundedness
of certain singular integral operators, especially on domains
in complex space.
\end{abstract}

\date{Revised July 15, 1995}

\maketitle

\section{Introduction} 

Function spaces play a significant role in harmonic analysis and
partial differential equations. The integral operators that form a
bridge between function spaces and partial differential equations are
the Calder\'on-Zygmund operators.  It is well known that
Calder\'on-Zygmund operators are bounded on the Lebesgue space
$L^p(\RR^n)$ for $1<p<\infty$. It is also known that the
Calder\'on-Zygmund operators are not bounded on $L^p(\RR^n)$ for any
$0<p\le 1$. It is natural to ask what are the substitutes for
$L^p(\RR^n)$ when $p$ is in this range; this circle of ideas has 
received considerable attention in harmonic analysis during the past
thirty years (see, for example, [COI], [CHR2], [COW1, 2], [FS], [KRA2]
and [STE], etc.).

We now understand that the best substitutes for the $L^p$ spaces are
the atomic Hardy spaces. Of course the holomorphic Hardy spaces on a
domain in $\CC^n$ and the real variable harmonic Hardy spaces on
$\RR^N$ are, at least on a formal level, quite different. It is
natural to wish to find a way to connect them. With this end in view,
the abstract Hardy spaces on a space of homogeneous type have been
introduced and studied by Coifman and Weiss [COW1, 2] and others.
For the case $0<p\le 1$, the atomic Hardy spaces $H^p(X)$ on a space
of homogeneous type were introduced in [COI] and [COW1, 2]. With
their definition of $H^p(X)$, it cannot be guaranteed that the
Calder\'on-Zygmund operators are bounded on $H^p(X)$ when
$p<1/2$---even when $X$ is the real line. 

One of the main purposes of this paper is to find a natural way to define
Hardy spaces $H^p(X)$ on a space of homogeneous type with $0<p<1$ in
such a way that singular integrals will be bounded on
all of these Hardy
spaces.  This will extend the work of Coifman and Weiss in [COW1],
[COW2].  In order to achieve this goal we shall have to address
several important ancillary issues:  how to define higher
order moment conditions on an arbitrary space of homogeneous type,
how to define analogues of smooth functions in that setting,
and how to define analogues of polynomials or other ``testing
functions''.

The second purpose of this paper is to establish duality theorems
for the Hardy spaces we have defined.
The study of singular integrals acting on $L^p$ spaces on a space of
homogeneous type has attracted many authors; one significant
recent study includes the adaptation of 
$T(1)$ or $T(b)$ theorems to $L^p(X)$ on
a space of homogeneous type $X$ (see [CW1, 2], [CHR1, 2] and
references therein.) 

The third purpose of this paper is 
to use some of those ideas to study the 
boundedness and compactness and other properties of some generalized Toeplitz operators (including commutators of a
singular integral operator and a multiplication operator) on function
spaces on a space of homogeneous type. In particular, we shall
extend the results in [CG], [CRW], [JAN] and [LI] from $\RR^n$ or $S^{2n-1}$
to much more general settings.
\medskip

The paper is organized as follows: In Section 2, we recall and prove
some preliminary results, and introduce some notation and
definitions. The atomic Hardy spaces $H^p(X)$ when $p<1$ are defined
and the duality theorem is stated and proved in Section 3.  In Section
4, we introduce several examples which fit our models. 
In Section 5, we prove boundedness of singular integral operators on 
$H^p$. In Section 6, we study the boundedness of some generalized
Toeplitz operators on 
$H^p$ or $L^p$.
The compactness of commutators on $L^p$ and $H^p(X)$ is 
studied in Section 7. 
Finally, in Section 8, we shall apply some of the theorems we proved 
in the previous sections
to prove analogous theorems for holomorphic Hardy space on some 
domains in $\CC^n$.
\medskip

\section{Preliminaries}

Let $X$ be a locally compact Hausdorff space.  A
 {\it homogeneous structure} on $X$ consists of a positive regular
Borel measure $\mu$ on $X$ and a family $\{B(x,r): x\in X, r > 0\}$
of basic open subsets of $X$ such that for some constants
$c>1$ and $K>1$ we have
\smallskip \\

(1) $x\in B(x, r)$ for all $x\in X$ and every $r>0$;

(2) If $x\in X$ and $0<r_1\le r_2$, then $B(x, r_1)\subset B(x, r_2)$;

(3) $0 < \mu(B(x, r)) < \infty$ for all $x\in X$ and all $r>0$;

(4) $X=\cup_{r>0} B(x,r)$ for some (and hence every) $x\in X$;

(5) $\mu(B(x,cr))\le K \mu(B(x,r))$ for all $x\in X$ and all $r>0$.

(6) If $B(x_1,r_1)\cap B(x_2, r_2)\ne \emptyset$ and $r_1\ge r_2$, then
$B(x_2, r_2)\subset B(x_1, cr_1)$;
\smallskip \\

We say that $X$ is a {\em space of homogeneous type} if $X$ is a locally
compact Hausdorff space having a {\it homogeneous structure}.  
Following Christ [CHR1], we assume from now on that
$$
 \mu(\{x\})=0   \leqno(2.1)
$$
for all $x\in X$.  

If $X$ is a space of homogeneous type, then 
one may define a quasi-distance on
 $X$ as follows: If $x, y\in X$, then
we let
$$
d(x,y)=\inf\{t: y\in B(x,t),\hbox{ and } x\in B(y,t)\} \leqno  (2.2)
$$
It is clear that $d(x,y)=d(y,x)$, and $d(x,y)=0$ if and
only if $x=y$.
From the so-called ``doubling property'' (5),
we have that $d(x,z)\le C(c, K) (d(x,y)+d(y,z))$. 
Therefore $d(\,\cdot\,,\,\cdot\,)$
is a quasi-metric on $X$.  
Coifman and Weiss refer to this quasi-metric
as the ``measure distance''.

For certain purposes, it is useful to choose a quasi-metric such that
the measure of a ball $B(x,r)$ associated to the quasi-metric is
comparable to a fixed power $r^\gamma$ of its radius. 
Follows Theorem 3 in [MS1],
, we have the following result:

\begin{lemma} \sl Let $(X, \mu)$ be a space of homogeneous type.
For any positive number $\gamma$, there is a quasi-metric 
$d_\gamma$ on $X$ such that if
$$
B_\gamma (x_0, r) = \{x \in X: d_\gamma (x, x_0) < r\} ,
$$
then
$$
\mu(B_\gamma (x_0, r)) \approx r^\gamma. \leqno(2.3)
$$
and with this quasimetric, we have
$$
\int_{X\setminus B(x_0, t)} \mu(B(x, t))^{-s}d\mu(x)
\le C_s \mu(B(x_0, t))^{-s+1}\leqno(2.4)
$$
for all $s>1$.
\end{lemma}

Now we define the maximum mean oscillation on balls with fixed radius
$r$ as follows:
$$
M(r,f)=\sup_{x\in X}\left \{{1\over \mu(B(x,r))}
             \int_{B(x,r)} |f-m_B(f)| \, d\mu \right \}.\leqno(2.5)
$$
[Here $m_B(f)$ is the mean value of $f$ on the ball $B$.]

\begin{definition}\rm
Let $(X,d,\mu)$ be a space of homogeneous type 
(the homogeneous structure is given by
the quasi-metric $d$.)
Let $f\in L^1_{\rm loc} (X)$. We say that $f\in BMO(X)$ if
$$
\|f\|_{BMO} \equiv \|f\|_* =\sup_{0<r<\infty} M(r, f)<\infty.
$$
We say that $f\in
VMO(X)$ if $f\in BMO(X)$ and
$$
\lim_{r\to 0^+} M(r, f)=0.
$$
\end{definition}

\noindent Now we may define the atomic $H^1$ space as follows.

\begin{definition}\rm
Let $(X,d,\mu)$ be a space of homogeneous type. Let $a\in L^\infty (X)$. 
We say $a$ is an {\em atom} (or a 1-atom)
if there is a ball $B$ such that $\hbox{supp} (a)\subset B$ and
$$
\mbox{(i)}\ \ |a(x)|\le 1/\mu(B) \quad ; \quad \mbox{(ii)}
              \ \ \int_B  a(x) \, d\mu=0.
$$
\end{definition}

\noindent Then
$$
H^1(X)= \Big \{u=\sum_{j=1}^{\infty} \lambda_j a_j: a_j 
    \hbox{ are atoms and }
\{\lambda_j\}_{j=1}^{\infty}\in \ell^1,\ \lambda_j\ge 0 \Big \}
\leqno(2.6)
$$
with norm
$$
\|u\|_{H^1}=\inf\Big\{\sum_{j=1}^{\infty} \lambda_j:
 u=\sum_{j=1}^{\infty} \lambda_j a_j\Big\}.
\leqno(2.7)
$$
The following result was proved in [COW2] and later in [MS2]:

\begin{theorem} \sl Let $X$ be a space of homogeneous type. Then
\vskip 7pt

(i) $ \bigl [H^1(X)\bigr ]^*=BMO(X)$;

(ii) $\bigl [VMO(X)\bigr ]^*=H^1(X)$.
\end{theorem}

The question of how to define atomic $H^p(X)$ spaces
when $0<p<1$ is more complex.  For
the case when $p$ is very closed to $1$, it was treated
by Coifman and Weiss [COW1,2] (also see [MS2]).
Let $a$ be a bounded function on $X$ with support 
in some ball $B=B(x_0, r)$.
We say that
$a$ is a $p$-{\em atom} in the sense of Coifman and Weiss if
$$
\mbox{(i)} \ \ |a(x)|\le \mu(B) ^{-1/p} \quad ; \quad
\mbox{(ii)} \ \ \int_{X} a(x) \, d\mu(x)=0 .
$$

\noindent The atoms of Coifman and Weiss are natural for values of $p$
that are close to $1$ (where the elementary mean value zero
property suffices for the purpose of studying singular integrals); when
the value of $p$ is small, then the definition of Hardy space requires
a higher order moment condition and is unworkable on an
arbitrary space of homogeneous type.

In order to define a more natural $p$-atom when $p$ is small,
we need to develop appropriate machinery
to produce the necessary moment conditions for $p$-atoms.
With a view to using Campanato-Morrey theory (see [KRA2]), 
one essential question is therefore
how to define polynomials
on a space of homogeneous type. Let us begin by taking a new look
at the problem on $\RR^n$,
where the traditional method is to exploit the canonical
coordinate system to define monomials and polynomials. 

We shall begin by replacing the usual polynomials by a family
of testing functions that have certain of the functorial properties
of polynomials.

\begin{definition} \rm Let $X$ be a space of homogeneous type. We say that
there is a family ${\cal P}$ of `testing polynomials' on $X$ if 
there is a positive integer-valued function $n(i, x)$ of 
the non-negative integer $i$
and the point $x\in X$ ($n(0,x)=1$) and families
of functions
$$
{\cal P}(x_0)=\{1\} \cup \{{\p}_{ij}(x, x_0): j=1,\cdots, 
        n(i, x_0), i=1,2,\cdots\} .
$$ 
[It is helpful to think of $i$ as the degree of the monomial ${\p}_{ij}$.]
We let ${\cal P}=\cup_{x_0\in X} {\cal P}(x_0)$.  In practice
we will assume that ${\cal P}$ contains the constant functions.

Let $k$ be a positive integer and $x_0\in X$.  We say that $p(x,
x_0)$ is a {\em polynomial} associated to
$x_0$, of {\em degree} not exceeding $k$, if there are 
constants $c_{ij},\ 0\le i\le k, 1\le j\le n(i,x_0)$ with
$$ 
p(x,x_0)=\sum_{i=0}^k \sum_{j=1}^{n(i, x_0)}
c_{ij} {\p}_{ij}(x,x_0), \leqno(2.8) 
$$ 
with each $p_{ij} \in {\cal P}(x_0).$  We will usually
consider constant functions to be polynomials.
\end{definition}

Let $(X, d, \mu)$ be a space of homogeneous type. We shall use the
notation $(X, d, \mu, {\cal P})$ to indicate the 
existence of the family of `testing polynomials'. 
\smallskip \\

\begin{example}  \sl
 Let $X=\RR^N$, $d(x,y)=|x-y|$, and
let $\mu$ be ordinary $N$-dimensional Lebesgue measure.
Then we let $n(i, x_0)= {N + i - 1 \choose i - 1}$ 
and $n(0,x)=1$ for every $x$ as usual.
Finally, $p_{ij}(x, x_0)=p_{ij}(x)$ is defined as follows:
 $p_{i1}(x)=x_1 ^i,\,  p_{i2}(x) =x_1^{i-1} x_2,\, 
\cdots,\,  p_{i\, n(i)}(x)=x_N^i$ (note here that we are simply
enumerating all the monomials of degree $i$).
\endpf
\end{example}

\begin{definition} \rm
Let $(X,d, \mu)$ be a space of homogeneous
type.  Let $1\le q\le \infty$, and $0<\alpha<\infty$.  
We shall say that a locally
integrable function $\phi$ belongs to $L(\alpha, q)(X)$ if
$$
\|\phi\|_{\alpha, q}
 =  \sup_B \left \{\mu(B)^{-\alpha}\Biggl [
\inf \int_{B} |\phi(x)-m_B(\phi)|^q \, {d\mu(x)\over \mu(B)}
\Biggr ]^{1/q} \right \}<\infty. \leqno (2.9)
$$
 
\noindent Here the supremum is taken over all balls $B$ in $X$.

We say that $\phi\in \hbox{Lip}_\alpha (X)$, $\alpha > 0,$ if
$$
\|f\|_{{\rm Lip}_{\alpha}}=\sup\Biggl \{{|f(x)-f(y)|
               \over d(x, y)^{\beta}} :
               x, y\in X, x\ne y \Biggr \}<\infty.
\leqno(2.10)
$$
\end{definition}

The following proposition is due to Mac\'as and Segovia [MS1].

\begin{proposition}\sl
 Let $(X, d,\mu)$ be a space of homogeneous type. Let
$1\le q\le \infty$ and $0<\alpha<\infty$. If  $\phi\in L(\alpha,q)$,
then there is a function $\psi(x)$ such that
\vskip 7pt

(i) $\phi(x)=\psi(x)$  for almost  all $x\in X$; \quad and
\vskip 6pt

(ii) $|\psi(x)-\phi(x)|\le C\|\phi\|_{\alpha,q} \mu(B)^{\alpha}$ where $B$ 
is any ball containing both $x$ and $y$.
\end{proposition}
 
\noindent {\bf Note:}  If the $\alpha$ is very large, then 
the space $L(\alpha, q)$ may contain only
constants.  This is a saturation theorem in the sense
of Favard.  For example, if $X=\RR^n$ and $d(x,y)=|x-y|$ 
and $\mu$ is the Lebesgue
measure, then $L(\alpha, q)=\CC$ if $\alpha>1/n$.
\smallskip

In what follows, we shall identify the functions $\phi$ and $\psi$
in the last proposition without further comment.

As to the relationship between $L(\alpha, q)$ and
$\hbox{Lip}_{\alpha}(X)$, we prove the following
proposition:

\begin{proposition}\sl Let $(X, d_\gamma,\mu)$ be a space of 
homogeneous type.
Let $1\le q\le \infty$ and $0<\alpha<\infty$.  Then
 $\phi \in L(\alpha, q)$  if and only if
$\phi \in \hbox{Lip}_{\gamma \alpha}(X)$.
\end{proposition}

\proof Suppose that $\phi \in \hbox{Lip}_{\gamma \alpha}(X)$.  
Let $B$ be any ball in $X$.  Then we have
\begin{eqnarray*}
\lefteqn{{1\over \mu(B)}\int_B |\phi(x)-m_B(\phi)| ^q \, d\mu(x)}\\
&\le &\|\phi\|_{{\rm Lip}_{\gamma \alpha }}^q \int_B
\Big({1\over \mu(B)}\int_B d_\gamma (x,y)^{ \gamma\alpha} \, d\mu(y)\Big)^q 
\, {d\mu(x)\over \mu(B)}\\
&\le& \|\phi\|_{{\rm Lip}_{\gamma\alpha}}^q\int_B
\Big({1\over \mu(B)}\int_B r^{\gamma \alpha } \, d\mu(y)\Big)^q \, {d\mu(x)
\over \mu(B)}\\
&\le&\|\phi\|_{{\rm Lip}_{\gamma \alpha}}^q  r^{\gamma \alpha q} \\
\end{eqnarray*}
Therefore
$$
\left({1\over \mu(B)}\int_B |\phi(x)-m_B(\phi)| ^q \, d\mu(x)\right)^{1/q}
\le \|\phi\|_{{\rm Lip}_{\gamma\alpha}} r^{\gamma \alpha}\le
C\|\phi\|_{{\rm Lip}_{\gamma\alpha}} \mu(B)^{\alpha}.
$$
Therefore $f\in L(\alpha, q)(X)$.

Conversely, let $\phi\in L(\alpha, q)$.  We will show that
$\phi \in {\rm Lip}_{\gamma\alpha}(X)$. By Proposition 2.8, we have
$$
|\phi(x)-\phi(y)|\le C\|\phi\|_{\alpha, q} \mu(B)^{\alpha}
$$
for any ball $B(x_0, r)$ containing $x$ and $y$. In particular, we have
$$
|\phi(x)-\phi(y)|\le C\|\phi\|_{\alpha, q} 
\mu(B(x, cd_\gamma (x,y))^{\alpha}
\le C^{\alpha}\|\phi\|_{\alpha, q} K^{\alpha} 
d_\gamma (x,y)^{ \gamma \alpha }
$$
This implies that $\phi \in {\rm Lip}_{\gamma \alpha}(X)$.
 Therefore the proof of the
proposition is complete.\endpf

Next we shall define smooth functions on $X$ by comparing
them with our testing functions above.

For $\gamma > 0$ fixed,
let $\beta_\gamma$ be the supremum of positive numbers $t_\gamma$ such that 
$$
|d_\gamma(x,z)-d_\gamma (y, z)|\le C(K,c) r^{1-t_\gamma}
 d_{\gamma}(x,y)^{t_\gamma}
$$
for all $x, y\in X$ with $d_\gamma(x, z)<r$ and $d_\gamma(y,z)<r$ 
and $r>0$. 
In order to
define smooth functions based on the family of testing polynomials, 
it is natural to
first assume that ${\cal P}\subset \hbox{Lip}_{\beta_\gamma}(X)$. 
Thus, we shall give the following definition.

\begin{definition} \rm
Let $(X,d_\gamma, \mu, {\cal P})$
be a space of homogeneous
type. Fix a number $1\le q\le \infty$. For $0<\alpha<
\gamma^{-1}\beta_\gamma $ 
and a nonnegative integer $k$,
we say that a locally integrable function $\phi$ belongs to 
$L(\alpha, k, q)$ if $B=B(x_0, r)$ and
$$
\|\phi\|_{\alpha, k, q}=
\sup_{B}\left \{ \Big [ {1\over \mu(B)^{\alpha+k/\gamma}}\inf_p \Big[
\int_{B} |\phi(x)-p(x,x_0)|^q \, {d\mu(x)\over \mu(B)}
\Big]^{1/q}\right \}<\infty . \leqno(2.11)
$$
Here the infimum is taken over all testing polynomials $p$ centered at
$x_0$ of degree not exceeding $k$, 
and the supremum is taken over all balls $B=B(x_0, r) \subseteq X$.
\end{definition}

\begin{proposition}\sl Let $(X, d_\gamma,\mu,{\cal P})$ 
be a space of homogeneous type. Let
$1\le q\le \infty$ and $0<\alpha<\infty$. Then 
$L(\alpha, 0, q)=L(\alpha, q)$.
\end{proposition}

\proof It is obvious that
$$
\|\phi\|_{\alpha,0, q} \le \|\phi\|_{\alpha,q} .
$$
The opposite inequality follows from the following standard fact.
$$
\Biggl (\int_{B(x_0, r)}|\phi-m_B(\phi)|^q\, {d\mu\over \mu(B)}\Biggr )^{1/q}
\le 2\inf\Biggl\{\Big(\int_B |\phi-c|^q \, {d\mu\over \mu(B)}\Big)^{1/q}:
 c\in \CC \Biggr\} 
$$
and the proof of the proposition.\endpf

From the definition of $L(\alpha, k, q)(X)$,
it is clear that $$L(\alpha, k, q_2)(X)\subset L(\alpha, k, q_1)(X)$$
if $q_2\ge q_1$.
We shall show in the next lemma that the converse is true as well.

\begin{lemma} \sl Let $(X, d_\gamma,\mu, {\cal P})$ be a space of
homogeneous type. Then we have $ L(\alpha, k, 1)(X)= L(\alpha, k,
\infty)(X)$ for all $0<\alpha< \gamma^{-1}\beta_\gamma$.  
\end{lemma}

Let us begin by recalling the following covering lemma. We will use
it in the proof of Lemma 2.11 and in later sections as well.

\begin{lemma}\sl
Let $(X, d_\gamma, \mu)$ be a space of
homogeneous type and $E$ be a compact subset of $X$. Suppose that there is
a family ${\cal F}$ of balls which cover $E$.  Then there is a
sequence $\{B(x_k, r_k)\}$ of disjoint balls such that $E$ is covered
by the family of balls $\{ B(x_k, 5c r_k)\}$. Moreover, there is a constant
$C$ such that for any point $x\in X$ there are at most $C$ balls
$B(x_k, 5c r_k)$ that contain $x$ (that is, the covering
$\{B(x_k, 5c r_k)\}$ has valence $C$).
\end{lemma}

\noindent  For details on the covering lemma, see [COW1], [COW2], or [CHR1].

\medskip

\noindent{\bf Proof of Lemma 2.11:} It suffices to prove that for any 
$f\in  L(\alpha, k, 1)(X)$ 
we have
$f\in  L(\alpha, k,\infty)$.  In fact we show that there is a function
$g\in {\rm Lip}_{\alpha \gamma}$ such that $g=f$ for a.e.\ $x\in X$
and for any ball $B(x_0, r)$
there is a testing polynomial of degree not exceeding $k$ so that
$$
|g(x)-p(x,x_0)|\le C\|f\|_{L(\alpha, k, 1)} \mu(B)^{\alpha+k/\gamma}
\leqno(2.12)
$$
for all $x\in B(x_0, r)$.

By Propositions 2.7 and 2.8, there is a function
$g\in\hbox{ Lip}_{\gamma \alpha}(X)$ such that $f(x)=g(x)$
for a.e.\ $x\in X$ and
$$
|g(x)-g(y)|\le C\|f\|_{\alpha, 0, q} 
 d_\gamma(x,y)^{\gamma\alpha},\quad x, y\in X.
$$
[Note here that we use only the result of 2.7 for the $L(\alpha,0,1)$
space---that is, only the case $k = 0.$]  Therefore, it suffices 
to prove that this $g$ satisfies (2.12).

Since $f\in L(\alpha, k, 1)$, we know that $g\in L(\alpha, k, 1)$.
For any ball $B(x_0,r)$
in $X$, there is a testing polynomial $p_1(x, x_0)$ of degree
less than or equal to $k$
such that
$$
\int_B |g(x)-p_1(x, x_0)| \, {d\mu\over \mu(B)}\le 
C\|f\|_{\alpha, k, 1}\mu(B)^{\alpha+k/\gamma} .
$$
Applying the facts that $g\in {\rm Lip}_{\alpha \gamma}(X)$
and $p_1(x,x_0)\in {\rm Lip}_{\beta_\gamma}$ (by assumption) and 
Proposition 2.7 again, we see that 
$$
|g(x)-p_1(x, x_0)-g(x_0)+p_1(x_0, x_0)|\le
 C\|f\|_{\alpha, k, 1}\mu(B(x_0, d_\gamma(x,x_0))^{\alpha+k/\gamma}
$$
for all $x\in B(x_0, r)$. 

Now we choose
$$
p(x, x_0)=p_1(x, x_0)+g(x_0)-p_1(x_0, x_0) .
$$
Then
$$
|g(x)-p(x, x_0)|\le C\mu (B)^{\alpha+k/\gamma}
$$
for $x\in B(x_0, r)$ and $r>0$. Here we have used
the covering lemma; the constant depends
only on $k$, and is independent of $x_0$.

We have completed the proof of (2.12), and thus of Lemma 2.11.
\endpf

\section{Hardy Spaces with $p<1$}

We begin with the definition of atom:

\begin{definition} \rm
Fix a positive number $n$.
Let $(X,d_n,\mu, {\cal P})$ be a space of homogeneous type ($d_n=d_\gamma$ 
when $\gamma=n$ in Lemma 2.1.)
Let $0<p\le 1$. Then a bounded  measurable function $a$
defined on $X$ is said to be
a $p$-{\em atom} if there is a ball $B=B(x_0, r)$ such that
\vskip 7pt 

(i) $\supp (a)\subset B$; \ \ (ii) $|a(x)|\le 1/\mu(B)^{1/p}$;

(iii) For any testing polynomial $p(x, x_0)$ of degree less than or equal 
to $[n(1/p-1)]$, we have
$$
\int_B a(x) p(x,x_0) \, d\mu(x)=0
$$
where $[x]$ is the greatest integer that does not exceed $x$.
\end{definition}

Now we are ready to define the real Hardy spaces
on $X$ when $p \leq 1$.

\begin{definition} \rm Let $(X,d_n,\mu, {\cal P})$ be a space of
homogeneous type.
Let $0<p\le 1$. A measurable function $f$ is said to belong
to $H^p(X)$ if there are a sequence
$\{\lambda_k\}\in \ell^p$ of non-negative numbers and
a sequence of $p$-atoms $\{a_k\}$ such that
$$
f=\sum_{k=1}^{\infty} \lambda_k a_k
$$
in the sense of distributions (when $p = 1$ we may actually
take the convergence to be in $L^1$ but when $p < 1$ we
must use the distribution topology). We let
$$
\|f\|_{H^p}=\inf \left \{\sum_{k=1}^{\infty} \lambda_k^p: f=\sum \lambda_k a_k \right \}
$$
\end{definition}

For simplicity, we shall use the following notation. For each $0<p\le 1$ and
positive number $n$, we let
$$
\alpha(p,n)=1/p-1-{1\over n}\bigl [n(1/p-1)\bigr ] .
$$
The main purpose of this section is to prove the following theorem.

\begin{theorem}\sl Fix a positive number $n$. 
Let $(X, d_n,\mu, {\cal P})$ be a space of
homogeneous type. Let $0<p\le 1$ be such that $n\alpha(p,n)\le \beta_n$. 
Then
$$
\bigl [ H^p(X) \bigr ]^*=
 L\left (\alpha(p,n), [n(1/p-1)], q \right )(X)
$$
for all $1\le q\le \infty.$
\end{theorem}

We shall break the proof of Theorem 3.3 up into several lemmas.  The
dimension parameter $n$ should be considered fixed once and for all.

\begin{lemma}\sl Let $(X,d_n,\mu, {\cal P})$ be a
space of homogeneous type. Let $0<p\le 1$ be such that 
$n\alpha(p,n)\le \beta_n$
 Then 
$$
L(\alpha(p,n), [n(1/p-1)], q)(X)\subset \bigl [ H^p(X) \bigr ]^*
$$ 
for all $1\le q\le \infty$.
\end{lemma}
\proof It suffices to prove the case $q=1$.
Let $f\in L(\alpha(p,n), [n(1/p-1)], 1)(X)$. 
Then, for any $u\in H^p(X)$ with
$$
u=\sum_{j=1}^{\infty} \lambda_j a_j,\quad \hbox{supp}(a_j)\subset B_j,
$$
we have for some $p(x,x_j)\in {\cal P}_{[n(1/p-1)]}(x_j)$ (the space 
of all testing polynomials, associated to $x_j$,
that have degree less than or equal $[n(1/p-1)]$ ) that
\begin{eqnarray*}
\lefteqn{\left|
\int_X f \ \sum_{j=1}^{k} \lambda_j a_j \, d\mu\right|}\\
&=&\left|
\sum_{j=1}^k\lambda_j \int_{B_j} (f(x)-p(x,x_j)  a_j(x) \, d\mu(x)
 \right|\\
&\le &
\sum_{j=1}^k\lambda_j \int_{B_j} |f(x)-p(x,x_j)| |a_j(x)| \, d\mu(x) \\
&\le &
\sum_{j=1}^k\lambda_j \mu(B_j)^{1-1/p}\int_{B_j} |f(x)-p(x,x_j)|
 { \, d\mu
\over \mu(B_j)} \\
&\le & C\|f\|_{\alpha(p,n), [n(1/p-1)], 1}
 \sum_{j=1}^k\lambda_j \mu(B_j)^{1-1/p}
 \mu(B_j)^{\alpha(p,n)+{1\over n} [n(1/p-1)]} \\
&\le & C\|f\|_{\alpha(p,n), [n(1/p-1)], 1}
\sum_{j=1}^k\lambda_j  \\
&\le & C\|f\|_{\alpha(p,n), [n(1/p-1)], 1} \|u\|_{H^p}.
\end{eqnarray*}
Therefore, letting $k\to \infty$, we have
$$
\left|\int_X f u \, d\mu\right|\le  C\|f\|_{L(\alpha(p,n), [n(1/p-1)], 1)} 
\|u\|_{H^p}
$$
This completes the proof of the lemma.\endpf

To prove the converse, we need the following notation and lemmas.

\begin{definition} \rm Let $(X, d,\mu, {\cal P})$ be a space 
of homogeneous type.
Let $a$ be a measurable function with support in some
 ball $B=B(x_0, r)$. We say
that $a$ is a $(p,2)$-{\em atom} if, instead of satisfying
the classical size condition 
$\|a\|_{\infty}\le \mu(B)^{1/p}$, the function $a$ instead satisfies
$\bigl (\int_B|a|^2 \, d\mu/\mu(B) \bigr )^{1/2}\le 1/\mu(B)^{1/p}$.
Let $H^{p, 2}(X)$ be the
atomic Hardy space created by replacing classical
$p$-atoms with $(p,2)$-atoms.
\end{definition}

Now we shall prove the following simple lemma.

\begin{lemma}\sl  Let $(X,d_n,\mu, {\cal P})$ be a space of 
homogeneous type.
Then $H^p(X)\subset H^{p,2}(X)$ and the embedding is continuous.
\end{lemma}
\proof  Let $u\in H^p(X)$. Then there
are a sequence of positive numbers $\{\lambda_k\}$ and a
 sequence of $p$-atoms
$\{a_k\}$ such that
$$
u=\sum_k \lambda_k a_k \quad ,\quad \|u\|_{H^p}^p\approx \sum_k\lambda_k^p .
$$
Since the $a_k$ are $p$-atoms, they must be $(p,2)$-atoms. 
Thus $u\in H^{p,2}(X)$
and
$$
\|u\|_{H^{p,2}}^p\le C\sum_k \lambda_k^p\le C\|u\|_{H^p}^p.
$$
Thus $H^p(X)\subset H^{p,2}(X)$ and the embedding is continuous. 
\endpf

Next we shall prove the following:

\begin{lemma}\sl  Let $(X,d_n,\mu, {\cal P})$ be a space 
of homogeneous type. Let $0<p\le 1$ be 
such that $n\alpha(p,n)\le \beta_n$.
 Then $ \bigl [ H^p(X)\bigr ]^*\subset L(\alpha(p,n), [n(1/p-1)], q)(X)$ for all
$1\le q\le \infty$.
\end{lemma}
\proof Let $\ell \in \bigl [H^p(X)\bigr ]^*$. We need to find 
a function $f\in L(\alpha(p,n), [n(1/p-1)], q)$
such that
$$
\ell(u)=\ell_{f}(u)=\int_X f u \, d\mu
$$
for all $u\in H^p(X)$.

Since  $ H^p(X)\subset H^{p,2}(X)$ and is dense, we may extend $\ell$
to be a  bounded linear functional on $H^{p,2}(X)$. More precisely,
each $H^{p,2}$ atom can be taken to be an $L^\infty$ function, so
elements of $[H^p]^*$ act naturally on these atoms [the indicated
extension of the functional from $H^p$ to $H^{p,2}$ is not valid just
by abstract functional analysis].   Alternatively, Coifman and Weiss
[COW2] have shown that $H^p$ and $H^{p,2}$ are equivalent spaces, so
the point may be taken as moot.  All of these matters are laid out in
detail in that source.

For each
ball $B=B(x_0, c_0)$, $\ell$ is a bounded linear functional on 
$L^2(B,p)$, the space
of all $L^2$ functions $u$ with supports on $B$ and satisfying
$$
\int_B u(x) p(x, x_0))\, d\mu(x)=0 .
$$
Here $p \in {\cal P}_{[n(1/p-1)]}(x_0)$.
Thus, by elementary Hilbert space theory,
there is an $f\in L^2(B)$ 
(depending on $x_0$) such that
$$
\ell(u)=\int_B f u \, d\mu,\quad \hbox{ for all } u\in L^2(B, p).
$$
Since $x_0\in X$  and $c_0>0$ are arbitrary,
there is a function $f\in L^2_{{\rm loc}}(X)$ such
$$
\ell(a)=\int_X f(x) a(x) \, d\mu(x)
$$
for all $p$-atoms. Therefore, we have
$$
\ell(u)=\int_X f(x) u(x) \, d\mu(x)
$$
for any $u$ a finite linear combination of $p$-atoms.

Therefore, by Lemma 2.14, it suffices to show that 
$$
f\in L(\alpha(p,n), [n(1/p-1)], 1)(X)
$$
and
$$  \|f\|_{\alpha(p,n), [n(1/p-1)], 1}
\le C\|\ell\|_{\bigl [H^p(X)\bigr ]^*}.
$$

For any $p$-atom $a$ with support in $B(x_0, r)$, we have 
$$
\ell (a)=\ell_f(a)=\int_X f(x) a(x)\, \, d\mu(x)
=\int_B \bigl ( f(x)-p(x, x_0) \bigr ) a(x) \, d\mu(x)
$$
for all $p(x,x_0)\in {\cal P}_{ [n(1/p-1)]}(x_0)$. 
Fix $x_0\in X$, and $r>0$. We have
\begin{eqnarray*}
\lefteqn{\inf\Big\{\int_{B(x_0,r)}|f(x)-p(x, x_0)| \, d\mu(x): 
p(x, x_0)\in {\cal P}_{[n(1/p-1)]}(x_0)\Big\}}\\
&=&\sup \Big\{\Big| \int_{B(x_0, r)} f(x) g(x) \, d\mu\Big|:
 \|g\|_{\infty}\le 1\Big\} ,
\end{eqnarray*}
where the infimum is taken over all testing polynomials associated to 
the ball $B(x_0, r)$, of degree not exceeding $[n(1/p]-1)]$; also
the supremum is taken on all $g\in {\bf P}(x_0, p)^{\perp}$, where 
\begin{eqnarray*}
\lefteqn{{\bf P}(x_0, p)^{\perp}=\Big\{ g\in L^{\infty}(B(x_0, r)):}\\
& &\int_{B(x_0, r)} p(x, x_0) g(x) \, d\mu(x)=0\ \ \hbox{ for all } p(x, x_0)\in
{\cal P}_{[n(1/p-1)]}(x_0)
\Big\}
\end{eqnarray*}
{\em We claim that} $g\in{\bf P}(x_0, p)^{\perp}$ and $\|g\|_{\infty}\le 1$ 
if and only if
$a(x)=\mu(B_r(x_0))^{-1/p} g(x)$ is a $p$-atom supported on $B_r(x_0)$. 
The proof of
this claim is obvious. Thus
\begin{eqnarray*}
\lefteqn{\mu(B)^{-1/p}\inf_p\{\int_{B}|f(x)-p(x,x_0)])| \, d\mu(x): 
p(x,x_0)\in 
{\cal P}_{[n(1/p-1)]}(x_0)\}}\\
&=&\sup\{\Big|
\int_{B}f(x)\, a(x) \, d\mu(x)\Big|: a \hbox{ is } p-\hbox{atom}\}\\
&\le & \|\ell\|_{\bigl [(H^p)\bigr ]^*}
\end{eqnarray*}
Therefore $f\in  L(\alpha(p,n), [n(1/p-1)], 1)(X)$.

Combining Lemmas 3.4 and 3.7, the proof of Theorem 3.3 is complete.\endpf

Before ending this section, let us make the following remark about the
number $\beta_\gamma$. It is easy to see that
$\beta_{\gamma}=\beta_{\gamma_0} \gamma_0/\gamma$
 for any $\gamma_0, \gamma>0$. Moreover, if $\gamma\ge \gamma_0$ and
$$
\gamma_0 \alpha-[\gamma_0 \alpha]\le \beta_{\gamma_0},\quad \alpha>0,
$$
then
$$
\gamma \alpha-[\gamma \alpha]\le \beta_{\gamma}
$$
since $\gamma [\gamma_0 \alpha]/\gamma_0\ge [\gamma \alpha]$. 
Therefore, for each
$0<p\le 1$, in order to guarantee that 
$\gamma \alpha(p,\gamma)\le \beta_\gamma$,
we may increase the value of $\gamma$ by a suitable amount.
\smallskip \\

\section{Some Examples}

In this section, we shall introduce several examples of
spaces of homogeneous type. In each of these examples,
one can see that  our definition for
the atomic Hardy space when $0<p\leq 1$ and our 
duality Theorem 3.3 are natural
generalizations of the classical theory.

\begin{example} \sl Let $X=\RR^N$, equipped with the Euclidean metric 
$d(x,y)=|x-y|$. It is clear that $d(x,y)\in \mbox{Lip}_1(\RR^{2N})$.
 Let $\, d\mu=dv$ be the Lebesgue measure on $\RR^N$. Then $(X, d, \mu)$
is a space of homogeneous type, and $\mu(B(x,r))\approx r^N $ 
for any $x\in \RR^N$ and $r > 0$. 

Now we define the family of testing polynomials.
Let $n(i, x)= {N + i - 1 \choose i - 1} =n(i)$. Let ${\p}_0(x)=1$, and 
$$
{\p}_{i1}(x,x_0)=x_1^i, \ {\p}_{i2}(x,x_0)= x_1^{i-1} x_2, \,\cdots,\, 
\ {\p}_{i n(i)}(x, x_0)=x_N^i
$$ 
for $i=1,2,\cdots$. Then
\medskip \\

(i) $L(\alpha-[\alpha], [\alpha])(\RR^N)=\Lambda_{N \alpha}(\RR^N)$,
the Zygmund class in $\RR^N$, for any $\alpha>0$;
\vskip 4pt

(ii) $a$ is a $p$-atom in our sense if and only if it is a $p$-atom in 
the classical sense,
i.e., there is ball $B=B(x_0, r)$ such that $\supp(a)\subset B$, 
$|a(x)|\le 1/[v(B)]^{1/p}$ and
$$
\int_{\RR^N} a(x) x^k dx=0,\qquad |k|=k_1+k_2+\cdots k_N\le [N(1/p-1)];   
$$
for all $  k=(k_1,\cdots, k_N)$ with $k_i\ge 0$ and  
$x^k=x_1^{k_1}\cdots x_N^{k_N}$.
\vskip 4pt

(iii) $\bigl [ H^p(\RR^n)\bigr ]^*=\Lambda_{n(1/p-1)}(\RR^n)$.
\end{example}

\noindent {\bf Verification of Example:} The proof of the results of
Example 2  follows from the main theorems of Krantz [KRA1] and Latter
[LAT].  \endpf

\bigskip
\begin{example} \sl In this example, we take $X$ to be
 the Heisenberg group $\HH_n$. 
We begin with some notation and definitions.

{\it The Heisenberg group}:
$\HH_n$ is the Lie group with underlying manifold 
$\RR^{2n}\times \RR $ and group operation ($z=x+iy,\  z'=x'+iy'$):
$$
b b'=(x, y,t) (x,' y',t')= (z,t)(z',t')=
( z+z', t+t'+ 2\hbox{Im} \langle z, z'\rangle).
$$ 
Haar measure on $\HH_n$ is, up to a constant multiple, 
$dv=dx dy dt$---Lebesgue 
measure on $\RR^{2n+1}$.  
For each $c>0$, we define the dilations on $\HH_n$ as follows:
$$
c\cdot (x, y,t)=(c\, x, c\, y, c^2 t).
$$
The norm of $(x, y,t)$ is defined as:
$$
\|(x, y,t) \|=\{ t^2 +(|x|^2+|y|^2)^2\}^{1/4}
$$
It is clear that $d((x,y,t), (x', y',t'))
\equiv \|(x-x', y-y',t-t')\|\in {\rm Lip}_{1/2}(\HH_n)$.
Thus $(\HH_n, d, dv)$ is a space of homogeneous type.

Moreover,
$$
\|c\cdot (x,y,t)\|=|c| \, \|(x, y,t)\|.
$$
For each $b\in \HH_n$ and $r>0$ we define the ball 
$$
B(b, r)=\{ h\in \HH_n: \|b h^{-1}\|<r\},\quad cB(b, r)=B(b, cr)
$$
Thus we have
$$
v(cB)= c^{2(n+1)} v(B),\quad c>0.
$$
Therefore
$$
v\bigl (B(g, r)\bigr ) \approx r^{2(n+1)},\quad g\in \HH_n,\ r>0.
$$

Let $g=(x_1,\cdots, x_n, y_1,\cdots, y_n, t)$. We say the $p(g)=p(x,y,t)$
is a polynomial of degree less than or equal $k$ if
$$
p(g)=\sum_{|\alpha|\le n} c_{\alpha} g^{\alpha}
$$
where $\alpha=( \alpha_1,\cdots, \alpha_n,\alpha_{n+1},\cdots,\alpha_{2n},
\alpha_0)$
with $\alpha_0$ and $\alpha_j, j\ge 1$  non-negative
integers, and
$$
|\alpha|=2\alpha_0+\sum_{j=1}^{2n}\alpha_j\le k.
$$
Such a function $p$ is an element of ${\cal P}_0$.
Then we define the family $\cal P$ of testing polynomials on $\HH_n$ as follows:
${\cal P}_k (x_0)={\cal P}_k(0)$ and ${\cal P}_k(0)$ is the set of all polynomials
 of degree not exceeding $k$. 
 Then ${\cal P}=\cup_{k=0}^{\infty} {\cal P}_k(0)$.

 Let $X=\HH_n$, the Heisenberg group with metric
and measure defined as the above. Moreover, we let 
the family ${\cal P}$ of `testing polynomials' be the set
${\cal P}$ of all polynomials
on $\HH_n$ with the definition of degree as above. Then
\vskip 6pt

(i) $L\bigl (\alpha-{1\over 2(n+1)}[2(n+1)\alpha],[2(n+1)\alpha]\bigr )(\HH_n)
= \tilde{\Lambda}_{(n+1) \alpha}(\HH_n)$ for any $\alpha>0$;
this is the non-isotropic Zygmund space;
\vskip 4pt

(ii) $a$ is a $p$-atom in our sense if and only if it is a 
$p$-atom in the ``classical'' sense,
i.e., there is a ball $B(h, r)$ in $\HH_n$ such that 
$\supp(a)\subset B(h, r)$, $|a(x)|\le 1/v(B(h, r))^{1/p}$ and
$$
\int_{\HH_n} a(x) g^\alpha dv=0,
$$
for all $ \alpha
\hbox{ and } |\alpha|=2\alpha_0+\sum_{j=1}^{2n} \alpha_j\le
 [2(n+1)(1/p-1)]$;
\vskip 4pt

(iii) $\bigl [H^p(\HH_n)\bigr ]^*=\tilde{\Lambda}_{(n+1)(1/p-1)}(\HH_n)$.
\end{example}

\proof Conclusion (i) follows from the main theorem 
in [KRA1] and (ii) is a corollary
of Theorem 3.3.\endpf

\begin{example} \sl
Let   
${\cal U}_n=\{(z, \xi)\in \CC^n\times \CC: \hbox{Im}\, \xi -|z|^2>0\}$ 
be the Siegel upper half space. Let $B_{n+1}$ 
be the unit ball in $\CC^{n+1}$. It is well known that ${\cal U}_n$ 
and $B_n$
are biholomorphically equivalent. 

The automorphism group of ${\cal U}_n$ decomposes
into three subgroups via the Iwasawa decomposition; the nilpotent
piece of the decomposition acts simply transitively on $\partial {\cal U}_n$.
Thus we may identify $\partial {\cal U}_n$ with that group, and
it is the Heisenberg group $\HH_n$.  By way of the biholomorphism 
(a generalized
Cayley transform), the boundary of the ball (less a point) may be
identified with the Heisenberg group.  [These matters are treated
in detail in [BCK].]

When the boundary of the ball is equipped with these 
``translations'', one has
the following results.  Let $X=\d B_{n+1}$. Then 
\medskip \\

(i)  $L(\alpha-{1\over 2(n+1)}[2(n+1)\alpha], [2(n+1)\alpha])
(\d B_{n+1})=\tilde{\Lambda}_{(n+1)\alpha}(\d B_{n+1})$ for all $\alpha>0$;
\smallskip \\

(ii) There are $p$-atoms on $\partial B$, which correspond to those
on the Heisenberg group;
\smallskip \\

(iii) $\bigl [ H^p(\d B_{n+1})\bigr ]^*=\tilde{\Lambda}_{(n+1)(1/p-1)}(\d B_{n+1})$
 \end{example}
 
\begin{example} \sl Let $X=\d \Omega$, where $\Omega$ is 
a strictly pseudoconvex domain in
$\CC^n$ with smooth boundary. 
Then $\partial \Omega$ has the structure of a space of homogeneous
type, and it is similar to that for the unit ball.
[See [KRA2], [KRA3] for details of this matter.]

Then results analogous to (i) - (iii) of the last example hold in
this context.
\end{example}
\smallskip

\begin{example} \sl Our Theorem 3.3 holds when $X$ is a
 {\em nilpotent} group discussed in [ROS] and [KRA1] by choosing
corresponding family of polynomials,
we omit the details here. 
\end{example}

\section{Singular Integrals on $H^p(X)$}

In this section, we shall study the boundedness of some
 singular integrals on the Hardy spaces $H^p(X)$. 
First we recall several definitions
from M. Christ [CHR2].

\begin{definition} \rm  
Let $(X,d, \mu)$ be a space of homogeneous type.
A {\it standard kernel} is a function $K:X\times X\setminus \{x=y\}
\to \CC$ such that there exist $\epsilon > 0$, and $0<C<\infty$
satisfying
$$
|K(x,y)|\le {C\over \mu(\lambda(x,y))}\quad \hbox{ for all distinct } x, y\in X ;
\leqno(5.1)
$$
here
$$
\lambda (x,y)=\mu(B(x, d(x,y))),
\leqno(5.2)
$$
and
$$
|K(x,y)-K(x',y)|
+|K(y,x)-K(y,x')|
\le \left({d(x,x')\over d(x,y)}\right)^\epsilon {C\over \lambda(x,y)}\leqno(5.3)
$$
whenever $ d(x,y)\ge c d(x,x')$.
\end{definition}

\begin{definition}\rm 
A continuous linear operator 
$T:\Lambda_\delta\to \Lambda_\delta'$
is said to be {\em associated} to $K$ if $K$ is locally integrable 
away from the diagonal and
$$
\langle Tf, g\rangle =\int \int K(x,y) f(y) g(x) \, d\mu(y) \, d\mu(x)
$$
for all $f, g\in \Lambda_\delta$ whose supports are separated 
by a positive distance.
Here $\Lambda_\delta'$ denotes the dual space of $\Lambda_\delta$.
\end{definition}

\begin{definition}\rm
 A singular integral operator $T$ is a continuous
linear operator from $\Lambda_\delta$ to $\Lambda_\delta'$ for some
 $\delta\in (0,\delta_0]$,
which is associated to a standard kernel.
\end{definition}

The following theorem is due to Coifman and Weiss [COW1]; 
the basic arguments appear in [CAZ] and [HAS].

\begin{theorem} \sl Any singular integral operator which 
is bounded on $L^2$ is also bounded on $L^p$ for all $p\in (1,\infty)$,
is weak type (1,1), and is bounded on $BMO$.
\end{theorem}

The celebrated $T(1)$ theorem of David and Journ\'e gives  
necessary and sufficient 
conditions
to test when a singular integral operator is bounded on $L^2$. 
Let $T$ be a singular
integral operator.  We say that $T$ is {\em weakly bounded}  on $L^2(X)$ if
$$
|\langle T\phi, \psi\rangle|\le C\mu(B(x_0, r))
$$
for all $\phi, \psi \in \Lambda_\delta(B(x_0,r))$ supported in $B$ and
$$
\|\phi\|_{L^2}\|\psi\|_{L^2}\le C\mu(B)
$$
for all $x_0\in X$ and $r>0$.
 It is easy to show that if $T$ is bounded on $L^2(X)$ then 
$T$ is weakly bounded on
$L^2$. Conversely, we have the following $T(1)$ theorem of David 
and Journ\'e.

\begin{theorem} \sl Any singular integral operator $T$ 
is bounded on $L^2$ if and only the following holds:
\smallskip \\

(i) $T$ is weakly bounded on $L^2(X)$;

(ii) $T(1) \in BMO(X)$ and $T^*(1)\in BMO(X)$.
\end{theorem}

It is known that a singular integral being bounded on $L^2$ does not imply that
it is  bounded on $L^p$ for $p\le 1$. The natural question is:

\begin{question} \sl Let $T$ be a singular integral operator 
which is bounded on
$L^2(X)$. Is $T$ bounded from $H^p$ to $L^p$ for all $0<p\le 1$ ?
\end{question}

The answer to this query is known to be ``yes'' when $p\le 1$ and 
sufficiently close to $1$. That
is the following result of [COW2].

\begin{theorem}\sl Any singular integral operator 
that is bounded on $L^2$ is also bounded from $H^p$ to  $L^p$  for
all $p\in (1-\epsilon_1,1]$
for some positive $\epsilon_1$ depending only on $X$ 
and the standard kernel $K$ associated to the operator.
\end{theorem}

From the proof of the above theorem, one can see that the 
$\epsilon_1$ depends on
the $\epsilon$ in the definition of a standard kernel $K$ 
since the atoms in $H^1$ have
$0$-order cancellation. In order to have a theorem like 
Theorem 5.6 when $p$ is small,
one can imagine that we need a condition on $K$ 
involving a higher order Lipschitz condition
that depends on $0<p<1$. Therefore we pose the following condition on a
standard kernel:

\begin{definition} \rm
Let $(X,d_n, \mu, {\P})$ be a space of homogeneous type.
A {\it $(p,k)$-standard kernel} is a function $K:X\times X\setminus \{x=y\}
\to \CC$ such that there exist $\epsilon>0$, a positive integer $k$,
 and $0<C<\infty$
satisfying (5.1), $\epsilon\ge 1/k$
and for fixed $x$ and $y$ 
there are two polynomials $p_1(x',y)$ and $p_2(x',y)$ in $x'$ with
 degree not exceeding $[n(1/p-1)] $ such that 
$$
|K(x',y)-p_1(x',y)|
+|K(y,x')-p_2(x',y)|
\le \left({d(x,x')\over d(x,y)}\right)^{\epsilon+{1\over k}[kn(1/p-1)]} 
{C\over \lambda(x,y)}   \eqno (5.4)
$$
whenever $d(x,y)\le c d(x,x')$. 
\end{definition}

Note that, when $[kn(1/p-1)]=0$, 
we take $p_1(x',y)=K(x,y)$ and $p_2(x', y)=K(y, x).$

The main purpose of this section is to prove the following theorem.

\begin{theorem}\sl Let $0<p\le 1$. Let $T$ be a singular 
integral operator with
$(p,k)$-standard kernel which 
is bounded on $L^2$.  Then $T$ is also bounded from $H^p$ to  $L^p$  for
all $0<p \le 1$.
\end{theorem}

\proof Let $T$ be a singular integral operator with a $(p,k)$-standard kernel $K(x,y)$.
To prove that $T:H^p(X)\to L^p(X)$ is bounded, 
from the definition of $H^p(X)$, it suffices
to prove that
$$
\|T(a)\|_{L^p}\le C_p
$$
for all $p$-atoms on $X$, where $C$ is a constant independent of $a$.
 Let $a$ be a $p$-atom with support in $B(x_0, r)$ and 
$$
\|a\|_{L^\infty}\le 1/\mu(B)^{1/p},\quad
\int_X  a(x) p(x,x_0) \, d\mu(x)=0
$$
for all $p(x, x_0)\in {\P}_{[n(1/p-1)]}(x_0)$. Thus
$$
T(a)(x)=\int_X K(x,y) a(y) \, d\mu(y)
$$
and
$$
\|T(a)\|_{L^2}=\left\|\int_B K(\cdot, y) a(y)\, d\mu(y)\right \|_{L^2}
\le C\| a \|_{L^2}=C\mu(B)^{(p-2)/(2p)} .
$$
Thus, for $1 \leq p < 2,$
\begin{eqnarray*}
\lefteqn{\left\|\int_B K(\cdot, y) a(y) \, d\mu(y)\right \|_{L^p}}\\
&\le& \left\|\chi_{c B}\int_B K(\cdot, y) a(y) \, d\mu(y)\right \|_{L^p}+
\left\|(1-\chi_{c B})\int_B K(\cdot, y) a(y) \, d\mu(y)\right \|_{L^p}\\
&\le& C\mu (c B)^{(2-p)/(2p)}\left\|\int_B K(\cdot, y) a(y)
 \, d\mu(y)\right \|_{L^2}\\
&&\quad +\left\|(1-\chi_{c B})\int_B K(\cdot, y) a(y) \, d\mu(y)\right \|_{L^p}\\
&\le& C( c, K)+
 \left\|(1-\chi_{c B})\int_B K(\cdot, y) a(y) \, d\mu(y)\right \|_{L^p}.
\end{eqnarray*}

Now we consider $x\in X\setminus c B$.  We choose
a polynomial $p_{x_0}(y, x) $ of degree less than or equal 
to $[n(1/p-1)]$ such that
\begin{eqnarray*}
|K(x, y)-p_{x_0}(y,x)|&\le&
 \left({d(y,x_0)\over d(y,x)}\right)^{\epsilon+[kn(1/p-1)]/k} {C\over \lambda(y,x)}\\
&\le& C
 \left({d(y,x_0)\over d(x,x_0)}\right)^{\epsilon+[kn(1/p-1)]/k} {C\over \lambda(x_0,x)}.
\end{eqnarray*}
Thus
\begin{eqnarray*}
\lefteqn{\left\|(1-\chi_{c B})\int_B K(\cdot, y) a(y) \, d\mu(y)\right \|_{L^p}^p}\\
&\le& \int_{X\setminus c B}
\Big|\int_B (K(x,y)-p_{x_0}(y, x)) a(y) \, d\mu(y)\Big|^p \, d\mu(x)\\
&\le& C\int_{X\setminus c B}
\Big[\int_B 
\left({d(y,x_0)\over d(x,x_0)}\right)^{\epsilon+[kn(1/p-1)]/k} 
{C \lambda(x_0,y)^{-1/p}\over \lambda(x_0,x)} \, d\mu(y)\Big ]^p \, d\mu(x)\\
&\le& C\int_{X\setminus c B}
\Big[
\left({ 1\over d(x,x_0)}\right)^{\epsilon+[kn(1/p-1)]/k} 
{C \mu(B)^{1-1/p+\epsilon/n+[kn(1/p-1)]/(kn)} \over \lambda(x_0,x)} \Big ]^p \, d\mu(x)\\
&\le& CC \mu(B)^{p-1+ p\epsilon/n+p[kn(1/p-1)]/(kn)}\int_{X\setminus cB}
\left({ 1\over d(x,x_0)}\right)^{p\epsilon+p[kn(1/p-1)]/k+np} 
 \, d\mu(x)\\
&\le& C \mu(B)^{p-1+ p\epsilon/n+p[kn(1/p-1)]/(kn)}\mu(B)^
{-p\epsilon /n -p[kn(1/p-1)]/(kn)-p+1} \\
&\le & C_p
\end{eqnarray*}
since $p\epsilon+ p[kn(1/p-1)]/k +np>n$. 

The proof of Theorem 5.8 is complete.
\endpf

It would be interesting to prove that the integral operator $T$ in
Theorem 5.8 is bounded on $H^p$ for $0<p\le$.  We shall leave the
matter for future study.

\section{Boundedness of generalized Toeplitz operators}

In this section, we shall study the boundedness  of
some generalized Toeplitz operators which become
from a family of singular integral operators (see (6.2)) in the context of
the Hardy spaces $H^p(X)$ and the Lebesgue spaces
$L^p(X)$ on a space of homogeneous type.

Let $f\in L^2(X,\mu)$. Then we define the multiplication operator 
$M_f$ with {\it symbol}
$f$ as follows: For any function $g$ on $X$, we let
$
M_f(g)(x)=f(x) \cdot g(x),\quad x\in X.$
Let $T_K$ be a singular integral operator with a standard kernel $K(x,y)$.
Then we denote by $C_f=[M_f,T_K]=M_f T_K-T_K M_f$, 
the commutator of $M_f$ and
$T_K$. In order to obtain better control over the singular integrals
defined in Section 5,
we let 
$$
K^\eta(x, y)=K(z, y) \ \hbox{ if } d(x,y)\ge \eta;\quad K^\eta(x,y)=0
\hbox{ if } d(x,y)<\eta
$$
and
$$
T^\eta (g)(x)=\int_X K^\eta(x,y) g(y) \, d\mu(y)
$$
and
$$
\tilde{T}(f)(x)=\sup_{0<\eta<<1} |T^\eta(f)(x)| .
$$
Then, by Theorem 12 in [CHR1], we have:
\medskip

$T$ is bounded on $L^2(X)$ $\iff$ $\tilde{T}$ is bounded on $L^p(X)$ for
$1<p<\infty$.
\medskip

 We
assume that the balls and measure $\mu$
satisfying the following condition: There is an $\epsilon_0>0$ such that
$$
c^{j\epsilon} \mu(B(x,t))\le \mu(B(x, c^j t)\le C(c, K)^j \mu(B(x,t)) \leqno(6.1)
$$
for all $t>0$ and $x\in X$. This condition is enough to guarantee 
(2.4) holds.

Let $T_{j,1}, T_{j, 2}$ $(j=1,\cdots, m$)
 are a finite sequence of C-Z type operators. We shall
consider a generalized Toeplitz operator:
$$
{\T}_b=\sum_{j=1}^m T_{j,1} M_b\, T_{j,2}. \leqno(6.2)
$$
We shall always assume that $T_{j,1}, T_{j,2}$ are bounded on
$L^2(X)$.

One of the main purposes of this section is to prove the following theorem.

\begin{theorem}\sl  Let $(X,d,\mu)$ be a space of homogeneous type.
Let $T_{j,i}$ are a sequence of C-Z operators which are bounded on
$L^2(X)$. If $g \in L^p(X)$ and ${\T}_1 (g)=0$, then for any
 $b\in BMO(X)$, we have ${\T}_b (g)\in L^p(X)$. Moreover, 
$$
\|{\T}_b(g)\|_{L^p(X)}\le C_p(\sum_{j=1}^m\|T_{j,1}\| )(\sum_{j=1}^m
\|T_{j,2}\|) \|g\|_{L^p} \|b\|_*
$$
for all $1<p<\infty$
\end{theorem}

By Theorem 5.4, we have immediatly the following corollary.

\begin{corollary}\sl  Let $(X,d,\mu)$ be a space of homogeneous type.
Let $T_{j,i}$ are a sequence of C-Z operators which are bounded on
$L^2(X)$ and ${\T}_1=0$. If $b\in BMO(X)$, then ${\T}_b$  is bounded on
$L^p(X)$ for all $1<p<\infty$.
\end{corollary}
\medskip

When $X$ is Heisenberg group, the above result was proved by 
L. Grafakos and  X. Li [GrL].

If we choose $m=2$ and $T_{1,1}=T_{2,2}=I$ and $T_{1,2}=T_{2,1}=T_K$, 
then ${\T}_b=C_b$. Thus

\begin{corollary}\sl  Let $(X,d,\mu)$ be a space of homogeneous type.
Let $T_K$ be a C-Z operator which is bounded on
$L^2(X)$. If $b\in BMO(X)$, then $[M_b,T_K]$  is bounded on
$L^p(X)$ for all $1<p<\infty$.
\end{corollary}

For each  $f\in L^1_{{\rm loc}}(X,\mu)$. 
Then we define
the sharp maximal function on $X$ as follows: For a.e.\ $x\in X$,
$$
f^{\#}(x)=\sup\Big \{{1\over
\mu(B)}\int_{B}|f(y)-f_{B}|\, d\mu(y): r>0 \Big \} .
 \leqno(6.3) 
$$
Also the $q$-maximal function of $f$ is this:
$$
 M_q(f)(x)=\sup\Big \{ \Big({1\over \mu(B)}\int_{B}|f(y)|^q
\, d\mu(y)\Big )^{1/q}:r>0 \Big \}.\leqno(6.4)
$$
Then we have  following three lemmas 

\begin{lemma}\sl  Let $(X, d,\mu)$ be a space of homogeneous type. Let
$f\in L^{p_0}(X,\mu)$ for some $1\le p_0<p$. 
Then $f\in L^p(X,\mu)$ if and only if $f^{\#}(x)\in L^p(X,\mu)$; and $f\in
L^p(X,\mu)$ if and only if $M(f)=M_1(f)\in L^p(X,\mu)$ for all $1<p<\infty.$
\end{lemma}

The proof of Lemma 6.4 can be found in M. Christ and R. Fefferman
[CHF]; Calder\'on [CAL] for $\Omega=R^n$; 
H. Aimar
and R. Maci\'as [AIM] for the Hardy-Littlwood maximal function on spaces of
homogeneous type. Also see J. O. Str\"omberg and A. Torchinski [STT] for sharp
maximal functions on spaces of homogeneous type.

\begin{lemma} \sl  
(a) $T$ is bounded on $L^2(X)$ $\iff$ $\tilde{T}$ is bounded on $L^p(X)$ for
$1<p<\infty$.

(b)  Let $1<p<\infty.$ Then $M_q(f)\in L^p(X,\mu)$
for all $1\le q<p$.

(c) If  $f\in BMO(X)$, then we have 
$
|f_{2^kB}-f_B|\le K\| f \|_* k
$
and
$$
\sup \Big \{ {1\over |B|}\int_B |f(y)-f_B|^p \, d\mu(y): B=B(x_0,r)\subset
X \Big\} < C_p \|f\|_*^p,
$$
where $B=B(x_0,r)\subset X$.
\end{lemma}

\medskip

Now we are ready to prove Theorem 6.1. 
\medskip

Let  $b\in BMO(X)$ have compact support, and $g\in L^p(X)$ with ${\T}_1(g)=0$.
Then ${\T}_b(g)\in L^{p_0}(X,\mu)$ for some
$1\le p_0<p$. By Lemma 6.4, it suffices to prove that ${\T}_b(g)^{\#}\in
L^p(X,\mu)$  
and $||{\T}_b(g)^{\#}||_p\le C_p||b||_*||g||_p$.

 Let $B=B(x,r)$ be any ball in $X$ and $cB=B(x, cr)$. 
 We let
$$
{\X}^1={\X}_{2B};\quad {\X}^2=1-{\X}_{2B}.
$$
Since ${\T}_1(g)=0$, and so ${\T}_{b_B}(g)=b_B {\T}_1(g)=0$. Thus
$$
{\T}_b (g)={\T}_{(b-b_B) {\X}^1} (g)+{\T}_{(b-b_B){\X}^2}(g)=g_1+g_2.
\leqno(6.5)
$$ 
Note that 
$$
g_1(y)={\T}_{(b-b_B){\X}^1}(g)(y)\\
=\sum_{j=1}^m T_{j,1}[(b-b_B) {\X}^1 T_{j,2}(g)](y).
$$
Thus for each $1<q<p$, we can choose $1<\gamma<\infty$ 
such that $q\gamma<p$.
\begin{eqnarray*}
\lefteqn{\Big(\int_{B}|g_{1}(y)
|^q d\mu\Big)^{1/q}}\\
&\le &\sum_{j=1}^m  C_{q,j} 
\Big( \int_{2B} |b-b_B|^q |T_{j,2}( g)(y)|^q d\mu(y)^{1/q}\\
&\le& \sum_{j=1}^m  C_{q,\gamma, j}
 \Big(\int_{B}|b-b_{B}|^{q
\gamma'}\, d\mu\Big)^{1/(q\gamma')}
\Big(\int_{2B}|T_{j,2}(g)|^{q\gamma}\, d\mu\Big)^{1/(q\gamma)}\\
&\le& \sum_{j=1}^m C_{q,\gamma, j}||b||_* \mu(B)^{1/(\gamma'q)}
M_{q\gamma}(T_{j,2}(g)(x) \mu(2B)^{1/q\gamma}\\
&\le&  ||b||_* \sum_{j=1}^m  C_{q,\gamma,j} M_{q\gamma}[T_{j,2}(g)](x)
\mu(B(x,r))^{1/(q\gamma)+1/(q\gamma')}\\
&\le &||b||_* \sum_{j=1}^m C_{q,\gamma,j} M_{q \gamma}[T_{j,2}(g)](x) 
\mu(B)^{1/q}.
\end{eqnarray*}
Therefore
$$
\int_{B}|g_{1}| \,{ d\mu \over \mu(B)}
\le \Big(\int_{B}|g_{1}|^q \, d\mu)^{1/q}\mu(B)\Big)^{1/q'-1}
\le \|b\|_*\sum_{j=1}^m  C_{q,\gamma,j} M_{q\gamma}[T_{j,2}(g)](x).
$$
Now we consider $g_{2}$. 
We shall prove the following lemma first.
\begin{lemma}\sl Let $T_K$ be a C-Z operator with a standard kernel $K$
such that $T_K$ is bounded on $L^2(X)$. Then for any $y\in B= B(x,r)$, 
we have
$$
|T_K(g{\X}^2)(y)-T_K(g{\X}^2)(x)|
\le C_{\epsilon,\epsilon_0} M(g)(x)
\leqno(6.6)
$$
and 
$$
|T_K[(b-b_B){\X}^2 g](y)-T_K((b-b_B){\X}^2 g](x)|
\le C_{\epsilon,\epsilon_0,\gamma}\|b\|_* M_{\gamma}(g)(x)
\leqno(6.7)
$$
where $ 1<\gamma<p$.
\end{lemma}
\proof The proof of (6.6) is similar and
 easier than the proof of (6.7). We shall present the proof of (6.7)
here.  Let $y\in B(x,r)$. Then
\begin{eqnarray*}
\lefteqn{|T_K[(b-b_B) g{\X}^2](y)-T_K[(b-b_B)g{\X}^2](x)|}\\
&=&\Big|\int_X (b-b_B) g(z){\X}^2(z)(K(y,z)-K(x,z)) \, d\mu(z)\Big|\\
&\le &\int_{X-2B}|b-b_B| |g(z)||K(y,z)-K(x,z)|\, d\mu(z)\\
&\le& \sum_{k=2}^{\infty} \int_{2^k B-2^{k-1}B} |b-b_B||g(z)||K(y,z)-K(x,z)|
\, d\mu(z)\\
&\le &C\sum_{k=2}^{\infty} \int_{2^k B}|b-b_B| |g(z)| \, d\mu(z)
\mu(B)^{\epsilon} \mu(B(x, 2^{k-1}r))^{-1-\epsilon}\\
&\le &C\sum_{k=2}^{\infty}\Big(\int_{2^k B}{|b-b_B|^{\gamma}\,
 d \mu \over \mu(2^k B)}\Big)^{1/\gamma'}
\Big (\int_{2^k B} {|g(z)|^{\gamma} \, d\mu(z)\over 
\mu(2^k B)}\Big)^{1/\gamma}
\mu(B)^{\epsilon} \mu(B(x, 2^{k-1}r))^{-\epsilon}\\
&\le &C\sum_{k=2}^{\infty} C_{\gamma'}\|b\|_* k M_{\gamma}(g)(x)
\mu(B)^{\epsilon} \mu(B(x, 2^{k-1}r))^{-\epsilon}\\
&\le & C_{\gamma}\|b\|_* M_{\gamma} (g)(x)
\sum_{k=2}^{\infty} \mu(B(x,r))^{\epsilon} k
 c^{-(k-2)\epsilon_0 \epsilon}
\mu(B(x,r))^{-\epsilon}\\
&\le & C_{\gamma} \|b\|_* M_{\gamma}(g)(x)
\sum_{k=2}^{\infty} k c^{-(k-2)\epsilon_0\epsilon}\\
&\le  & CC_{\epsilon,\epsilon_0,\gamma } \|b\|_* M_{\gamma} (g)(x),
\end{eqnarray*}
and the proof of lemma is complete.\epf

Since
$$
{1\over \mu(B)}\int_{B} |g_2(y)-(g_2)_{B}| \, d\mu(y)
\le {2\over \mu(B)}\int_B |g_2(y)- g_2(x)| d\mu(y).
$$
and
$$
g_2(y)=\sum_{j=1}^m T_{j,1}[(b-b_B){\X}^2 T_{j,2}(g)](y).
$$
 For each $y\in B(x,r)$, we have
$$
|g_2(y)- g_2(x)|
\le \sum_{j=1}^m  C_{\epsilon, \epsilon_0,\gamma,j}\|b\|_*
M_{\gamma}(T_{j,2}(g))(x)
$$
Thus
$$
{1\over \mu(B)}
\int_{B}|g_2(y)-(g_2)_B| d\mu(y)
\le \sum_{j=1}^m  C_{\epsilon, \epsilon_0,\gamma,j}\|b\|_* 
M_{\gamma}(T_{j,2}(g))(x)
$$
Since $T_{j,i}$ are bounded on $L^p(X)$ and
$$
{\T}_b(g)^{\#}(x)=g_1^{\#}(x)+ g_2^{\#}(x)
$$
Combining this and the estimation of $g_1^{\#}(x)$ and $g_2^{\#}(x)$,
we have
$$
\|{\T}_b(g)^{\#}\|_{ L^p}
\le C_p \sum_{j=1}^m \|T_{j,1}\|_{\infty} \|T_{j,2}\|_{\infty}\|b\|_*
\|g\|_{L^p(X)}.
$$
 By Lemmas 6.4, we have
${\T}_b(g)\in L^p(X,\mu)$ for all $1<p<\infty$. So the proof of Theorem 6.1 is
complete.
\epf

As a direct consequence of Theorem 6.1, we have the following theorem
in [CG] and [GL] for $X$ is $\RR^n$ or some homogeneous group.

\begin{theorem}\sl  Let $(X,d,\mu)$ be a space of homogeneous type.
If $f\in L^p(X)$ such that ${\T}_1(f)=0$, then the linear operator
$B_f(g)=\sum_{j=1}^m ( T_{j,1}^*(g), T_{j,2}(f))
$ is bounded from $ L^q(X)\to H^1(X)$,
where $p, q>1$ and $1/p+1/q=1$ and $(\cdot, \cdot)$
 denotes the inner product in $\CC^n$.  
\end{theorem}

\proof Since

\begin{eqnarray*}
\langle {\T}_b (f), g\rangle  
&=&\sum_{j=1}^m \langle T_{j,1}[b T_{j,2}(f)], g\rangle\\
&=&\sum_{j=1}^m \langle [b T_{j,2}(f)], T^*_{j,1}(g)\rangle\\
&=&\sum_{j=1}^m \langle b,\overline{ T_{j,2}(f)} T^*_{j,1}(g)\rangle\\
&=& \langle b,\sum_{j=1}^m \overline{ T_{j,2}(f)} T^*_{j,1}(g)\rangle\\
&=&\langle b, B_f(g)\rangle
\end{eqnarray*}
for all $b\in BMO(X)$ with compact support, and by Theorem 6.1, we have
$$
|\langle b, B_f(g)\rangle|\le C_{p,q}(\sum_{j=1}^n\|T_{j,1}\|)
(\sum_{j=1}^n\|T_{j,2}\|)\|b\|_* \|f\|_{L^p}\|g\|_{L^q}
$$
Since VMO functions with compact support is in $VMO(X)$ and Theorem 2.4,
we have $B_f(g)\in H^1(X)$, and the proof is complete.\epf

\medskip

The second purpose of this section is to  prove the following theorem.

\begin{theorem}\sl  Let $(X,d,\mu)$ be a space of homogeneous type. 
Let $T_K$ be 
a singular integral operator 
with a standard kernel $K(\cdot,\cdot)$. Suppose that $T_K$ is
bounded on $L^2(X)$ and that $f\in BMO(X)$. Then we have that
$C_f$ is bounded from $H^1(X)$ to $L_{\rm loc}^1(X)$.
\end{theorem}

\proof  From the definition of $H^1(X)$, it suffices to prove that
$$
\|C_f(a)\|_{L^1(X_0)}\le C\mu(X_0)\|f\|_*\leqno(6.8)
$$
for all atoms $a$ and any compact subset $X_0$ in $X$.
 Let $a$ be an atom with support $B=B(x,r)$.  This can be done
by using a standard computation, we omit the details here.
\endpf

In order to state and prove our next theorem, we reqire that our testing
polynomials $\cal P$ be closed under the product operation, 
more precisely, if $p(x_0, x)$
and $q(x_0, x)$ are polynomials of degree at most $n$ and $m$, then $pq$
is a polynomial of degree at most $n+m$.  This is analogous, but not
identical to, requiring that our testing polynomials form a graded
ring. Along this direction, we have 
 
\begin{theorem}\sl Let $0<p\le 1$ and
let $(X,d,\mu,{\P})$ be a space of homogeneous type. 
Further suppose that $\cal P$ is closed under
product operation.  Let $T_K$ be 
a singular integral operator 
with a $p$-standard kernel $K(\cdot,\cdot)$. If $T_K$ is
bounded on $L^2(X)$ and if
  $f\in L(\alpha(p,n),[n(1/p-1)],\infty)$, then we have that
$C_f$ is bounded from $H^p(X)$ to $L^p(X)$.
\end{theorem}

\proof  By Theorem 5.8, we have $\|T_K(g)\|_{L^p}\le C\|g\|_{H^p}$. Since
$$f\in L(\alpha(p,n),[n(1/p-1)],\infty),$$ we have $f\in L^\infty(X)$. Thus
$M_fT_K(g)\in L^p(X)$. Therefore the proof of Theorem 6.7 is reduced
to proving the following theorem.

\begin{theorem}\sl Let $0<p\le 1$ and
let $(X,d,\mu,{\P})$ be a space of homogeneous type. Let $T_K$ be 
a singular integral operator 
with a $p$-standard kernel $K(\cdot,\cdot)$. If $T_K$ is
bounded on $L^2(X)$, and if
  $f\in L(\alpha(p,n),[n(1/p-1)],\infty)$ then we have that
$T_K M_f$ is bounded from $H^p(X)$ to $L^p(X)$.
\end{theorem}

\proof By Theorem 5.8, it suffices to prove 
that $K(x, y) f(y)$ is a $(p,k)$-standard kernel. This is a 
direct computation
since $\cal P$ is closed under the product operation and
 $f\in L(\alpha(p, n), [n(1/p-1)], \infty)$. We omit the details.
\endpf

\section{Compactness of Commutators}

 In this section, we shall study the compactness of the commutator of
a singular integral operator and a multiplication operator on $L^p(X)$ with
$p>1$ on a space of homogeneous type. 

In order to study the compactness of $C_f$, we assume that $K\in C(X\times
X \setminus \{(x,x):x\in X\})$. We also assume that the measure 
$\mu$ satisfies the
following condition: There is a positive constant 
$\epsilon_0<\epsilon/2$ such that 
$$
|\mu(B(x,r)-\mu(y,r)|\le C \mu(B(x,d(x,y))^{\epsilon_0}  \leqno(7.1)
$$
for some $\epsilon_0>0$, all $x, y\in X$ and $d(x,y)\le r<1$.

The main purpose of this section is to prove the following theorem.

\begin{theorem}\sl  Let 
  $(X,d,\mu)$ be a space of homogeneous type satisfying (7.1). Let $T_K$ be 
a singular integral operator 
with a standard kernel $K(\cdot,\cdot)$. If $T_K$ is
bounded on $L^2(X)$, and
if $f\in VMO(X)$ then $C_f:L^p(X,\mu)\to L^p(X,\mu) $ 
is compact for all $1<p<\infty.$
\end{theorem}

Let $UC(X)$ denote the set of all uniformly continuous functions on $X$;
$BUC(X)$ denotes the subspace of all bounded uniformly continuous functions.
We begin by stating the following lemmas (their proofs are given later.) 

\begin{lemma}\sl Let $X$ be a space of homogeneous type. Let $f\in
VMO(X,\mu)$. Then for any $\eta>0$, there is a function
$f_{\eta}\in BUC(X)$
such that 
$$ 
||f_{\eta}-f||_* < \eta.
\leqno(7.2)
$$
Also there is an $\epsilon_1$, $0<\epsilon_1\le\epsilon_0/2<1$, such that
$$
|f_{\eta}(x)-f_{\eta}(y)|\le C_{\eta}\mu(B(x,r(x,y))^{\epsilon_1}.
\leqno(7.3) $$
\end{lemma}

For each $0<\eta<<1$, we let $K^{\eta}(x,y)$ be a continuous extension of
$K(x,y)$ from $X\times X-\{(x,y): d(x,y)<\eta\}$ to $X\times X$ such that
\begin{eqnarray*}
K^{\eta}(x,y) & = & K(x,y) \qquad \mbox{if} \ \  d(x,y)\ge \eta ;\\
|K^{\eta}(x,y)| & \le &
  C\mu(B(x,\eta))^{-1} \qquad \mbox{if} \ \ d(x,y)<\eta ; \\
 K^{\eta}(x,y) & = & 0 \qquad \mbox{if} \ \ d(x,y)\le \eta/c .
\end{eqnarray*}
We shall let $T_{K^\eta}$ denote the integral operator associated to $K_\eta$  
and we let $C^{\eta}_f=[M_f,T_{K^\eta}]$. Then we have the following Lemma.

\begin{lemma}\sl Let $f\in BUC(X)$ satisfy (7.2). Let  $T_K$ be a singular
integral operator with a standard kernel $K$ satisfying (7.1)
which is bounded on $L^2(X,\mu)$.
 Then $||C_f -C^{\eta}_f||_\infty\to 0$, as an operator on $L^p(X,\mu)$,
when $\eta\to 0$.
\end{lemma}

We shall postpone the proof of Lemmas 7.2 and 7.3. Let us first prove 
Theorem 7.1.
\medskip

\proof Applying Theorem 6.1, we have
$$
||C_f(g)-C_{f_{\eta}}(g)||_p\le C_p\|g\|_{L^p} ||f-f_{\eta}||_*<\eta.
\leqno(7.5)
$$
Therefore, in order to prove that $C_f$ is compact on $L^p(X,\mu)$, 
it suffices to prove that $C_{f_{\eta}}$ is compact on $L^p(X,\mu)$.

So, for convenience, we may assume that $ f\in BUC(X)$ satisfies 
(7.2) for some
constant $C$ and $\epsilon_0$. By Lemma 7.2, it suffices to prove 
that $C^{\eta}_f$ is compact on $L^p(X,\mu)$.

Since $K(x,y)\in C(X\times X \setminus \{(x,x): x\in X\})$, 
we have for each $g\in L^p(X)$ that $C^{\eta}_f(g)\in C(X)$. 
Moreover, for any $x, y\in X$ with $d(x,y)<1$, 
we have
\begin{eqnarray*}
\lefteqn{ C_f^{\eta}(g)(x)-C_f^{\eta}(g)(y)}\\
&=&f(x) \int_X K^{\eta}(x,z) g(z) \, d\mu(z)-\int_X K^{\eta}(x,z)f(z)g(z) \, d\mu(z)\\
&&\quad -f(y)\int_X K^{\eta}(y,z) g(z) \, d\mu(z)+\int_X K^{\eta}(y,z) f(z) g(z) \, d\mu(z)\\
&=&(f(x)-f(y))\int_X K^{\eta}(x,z) g(z) \, d\mu(z)\\
&&\quad +f(y) \int_X (K^{\eta}(x,z)-K^{\eta}(y,z)) g(z) \, d\mu(z)\\
&&\quad + \int_X (K^{\eta}(y,z)-K^{\eta}(x,z))f(z) g(z) \, d\mu(z)\\
&=&(f(x)-f(y)) \int_X K^{\eta}(x,z)g(z) \, d\mu(z)\\
&&\quad +  \int_X (K^{\eta}(y,z)-K^{\eta}(x,z))(f(z)-f(y)) g(z)\, d\mu(z)\\
&=& I_1(x,y)+ I_2(x,y).
\end{eqnarray*}
And
\begin{eqnarray*}
|I_1(x,y)|
&=&|(f(x)-f(y))|(|\int_X K^{\eta}(x,z)g(z) \, d\mu(z)|\\
&\le& C |f(x)-f(y)|
\{\int_{X}|K^{\eta}(x,z)|^{p'}\, d\mu(z)\}^{1/p'}||g||_p\\
&\le& C_{\eta,p'}|f(x)-f(y)| ||g||_p. 
\end{eqnarray*}

Notice that $f$ is bounded. If we let
$r=d(x,y)<(1/c^2)\eta$, then 
\begin{eqnarray*} 
|I_2(x,y)| 
&\le &\int_{X-B(x,\eta/c)}|K(y,z)-K(x,z)|| (f(z)-f(y)||g(z)|
\, d\mu(z)\\
&\le &\|f\|_{L^\infty}\int_{X-B(x,\eta/2)}C\mu(B(x, r))^{\epsilon}
\mu(B(x,d(x,z))^{-1-\epsilon} |g(z)|
\, d\mu(z)\\
&\le & C\|f\|_{L^\infty} \mu(B(x,d(x,y))^{\epsilon}
\biggl \{\int_{X-B(x,\eta/2)}\mu(B(x,r(x,z))^{-p'(1+\epsilon)}\, d\mu(z) \biggr \}^{1/p'}||g||_{L^p}
\\
&\le&  C_{\eta, p'} \mu(B(x,d(x,y))^\epsilon ||f||_{L^\infty} \|g\|_{L^p}.
\end{eqnarray*}
Therefore $\bigl \{C_f^{\eta}(U)\bigr \}$ is an
 equicontinuous family.  Here $U$ is the unit ball in
$L^p(X,\mu)$. Therefore the Ascoli/Arzela Theorem shows that 
$\bigl \{C_f^{\eta}\bigr \}$ is compact on
$L^p(X,\mu)$. This completes the proof of Theorem 7.1. 
\endpf

Next we return to the proofs of Lemmas 7.2 and 7.3.
\medskip

\def\PV{\hbox{PV}}

\noindent{\bf Proof of Lemma 7.3}
\medskip

Let $g\in L^p(X)$. Then for each $x, y\in X$, we have
\begin{eqnarray*}
\lefteqn{C_f(g)(x)-C^{\eta}_f(g)(x)}\\
&=&f(x)  \int_{B(x,\eta)} K(x,z) g(z) \, d\mu(z)-\int_{B(x,\eta)}
K(x,z)f(z)g(z) \, d\mu(z)\\
&&\quad -f(x)\int_{B(x,\eta)-B(x,\eta/c)} K^{\eta}(,z) g(z)
\, d\mu(z)+\int_{B(x,\eta)-B(x,\eta/c)}K^{\eta}(x,z) f(z) g(z) \, d\mu(z)\\
&=&-\int_{B(x,\eta)} K(x,z)(f(z)-f(x) g(z) \, d\mu(z)\\
&&\quad +\int_{B(x,\eta)-B(x,\eta/c)}
(K^{\eta}(x,z)(f(z)-f(x)) g(z) \, d\mu(z)\\
&=&I_1(x)+I_2(x).
\end{eqnarray*}
We first estimate $I_2(x)$, and we may let $\eta\le 1/c^2$.  Then
\begin{eqnarray*}
|I_2(x)|
&\le& \max\{|f(z)-f(x)|:z\in B(x,\eta)\}|\times\\
&&\quad\int_{B(x,\eta)-B(x,\eta/c^2)} 
C\mu(B(x,\eta/c^2))^{-1}
|g(z)| \, d\mu(z)|\\
 &\le&
C\max\{|f(z)-f(x)|:z\in B(x,\eta)\} M(g)(x). 
\end{eqnarray*}
 
\noindent Note that
\begin{eqnarray*} I_1(x) 
&=&|\PV\int_{B(x,\eta)}K(x,z) (f(z)-f(x))g(z) \, d\mu(z)|\\
&\le& \int_{B(x,\eta)} C\mu(B(x,d(x,z))^{-1} \mu(B(x,d(x,z))^{\epsilon_0}
|g(z)| \, d\mu(z)\\
&\le& C\mu(B(x,\eta))^{\epsilon_0/2}\int_{B(x,\eta)}
\mu(B(x,d(x,z))^{-1+\epsilon_0/2}|g(z)|\, d\mu(z)\\
\end{eqnarray*}
Applying Schur's Lemma, we have
$$
\int_{X}\left \{\int_{B(x,\eta)}\mu(B(x,r(x,z))^{-1+\epsilon_0/2} |g(z)|\, d\mu(z) \right \}^p
\, d\mu (x)\le C_{\epsilon_0}||g||_p^p.
$$
for all $\eta\le 1$.

Therefore we have
$$
||C_f(g)-C_f^{\eta}(g)||_p^p \le C^p \mu(B(x,\eta))^{\epsilon_0 p/4}||g||_p^p.
$$
for all $g\in L^p(X,\mu)$.

This implies that
$$
||C_f-C_f^{\eta}||\to 0,\quad \hbox{ as an operator in }L^p(X,\mu),\quad\hbox{
when } \eta\to 0.
$$
The proof of Lemma 7.3 is complete.\endpf

Finally, we prove Lemma 7.2.
In order to achieve this end, let us first prove the following Lemma.

\begin{lemma}\sl  Let $f\in BMO(X,\mu)$. Then we have
$$
|f_{B(x,r)}|\le C||f||_* \log (C/\mu(B(x,r)),\quad r\le 1;
\leqno(7.6)
$$
and 
$$
{1\over \mu(B(x,r))}\int_{B(x,r)-B(y,r)} |f| \, d\mu(z)\le C_{r}||f||_* \mu(B(x,
d(x,y))^{\epsilon_0}
\leqno(7.7)
$$
for all $x, y\in X$ and $ cr(x,y)\le r\le 1$.
\end{lemma}

\proof Let us prove (7.7) first. By hypothesis (7.3), we have
\begin{eqnarray*}
\lefteqn{\int_{B(x,r)-B(y,r)} |f| \, d\mu(z)}\\
&\le& (\int_{B(x,r)} |f(z)|^2 \, d\mu(z))^{1/2} (\int_{B(x,r)-B(y,r)} 1
\, d\mu(z))^{1/2}\\
&\le& C_r C \mu(B(x, d(x,y))^{\epsilon_0/2}.
\end{eqnarray*}
This completes the proof of (7.7). And the proof of (7.6) can be found in
[KRL2].\endpf

\noindent{\bf Proof of Lemma 7.2.}
\smallskip

For any $\eta>0$, since $f\in VMO(X,\mu)$, there is a $\delta(\eta)>0$ 
such that
$$
{1\over \mu(B(x,r))}\int_{B(x,r)}|f(y)-f_{B(x,r)}|\, d\mu(y)< \epsilon/C^2
$$
for all $r<c\delta(\eta)$.
Let
$$
f_{\eta}(x) =f_{B(x,\delta(\eta))}.
$$
Then we have
\begin{eqnarray*}
\lefteqn{|f_{\eta}(x)-f_{\eta}(y)|}\\
&\le& |{1\over \mu(B(x,\delta))}\int_{B(x,\delta)}
f(z)-f_{B(y,\delta)}\, d\mu(z)|\\
& &\quad +\left|{1\over B(x,\delta)}-{1\over
\mu(B(y,\delta)}\right|\int_{B(y,\delta)} f(z)\, d\mu(z) |\\
&=&I_1(x,y)+I_2(x,y) .
\end{eqnarray*}
Also, by Lemma 7.4, we have
$$
I_1(x,y)\le C_{\delta} \mu(B(x,d(x,y))^{\epsilon_0/2}.
$$
and 
\begin{eqnarray*}
I_2(x,y)&\le& |\mu(B(x,\delta)-\mu(B(y,\delta)|||f||_{*}{|\log \mu(B(y,
\delta)|\over \mu(B(x,\delta)}\\
&\le& C_{\delta} \mu(B(x, d(x,y))^{\epsilon_0/2}
\end{eqnarray*}
for all $x, y\in X$.
Therefore $f_{\eta}$ satisfies (7.3).

Next we show that (7.2) holds.

Let $B=B(x,r)$ be any ball on $X$ with $r>\delta$, let $\{ B(x,\delta): x\in
B\}$ be an open cover of $B$. By Theorem 1.2 in [COW2], there are 
pairwise disjoint
balls $B(x_j,\delta), j=1,2,\cdots, N$ such that $B\subset \cup_{j=1}^N
B(x_j,c\delta)$.  It is clear that $\cup_{j=1}^N
B(x_j,c\delta) \subset B(x,c_0cr)$, and $c_0$ is independent of $\delta$. 
Therefore we have
\begin{eqnarray*}
\lefteqn{{1\over \mu(B)}\int_B |f(y)-f_{\eta}(y)-f_{B}+{f_{\eta}}_B| \, d\mu(y) }\\
&\le& 2{1\over \mu(B)}\int_B |f(y)-f_{\eta}(y)| \, d\mu(y)\\
&\le& {2\over |B|}\sum_{j=1}^N \int_{B(x_j,c\delta)} |f(y)-f_{\eta}(y)|\, d\mu(y)\\
&\le& {1\over |B|}\sum_{j=1}^N \Big\{\int_{B(x_j,c\delta)}
|f(y)-f_{B(x_j,c\delta)}(y)|\, d\mu(y)\\
&&+\int_{B(x_j,c\delta)}|f_{B(x_j,c\delta)}-f_{\eta}(y)|\, d\mu(y)\Big\}\\
&\le& {2\over |B|}\sum_{j=1}^N (\epsilon/C)||f||_* \mu(B(x_j, c_0 c\delta))\\
&& +{2\over |B|}\sum_{j=1}^N \int_{B(x_j,c\delta)}{1\over
\mu(B(y,\delta))}\int_{B(y,\delta)}
 |f(z)-f_{B(x_j,c\delta)}|\, d\mu(z)\, d\mu(y)\\
&\le& (\epsilon/C^2)||f||_*{2\over |B|}\mu(B(x, c_0cr))
 +{2\over |B|}\sum_{j=1}^N \mu(B(x_j,c\delta))||f||_* C(\epsilon/C^2) \\
&\le& C ||f||_* \epsilon
\end{eqnarray*}
where $C$ is a constant independent of $\delta$ and so $\epsilon$.
Therefore the proof of Lemma 7.2 is complete. \endpf

\section{The Case that 
$X$ is Either $\RR^N$ or the Boundary of a Domain in $\CC^n$}

In this section, we shall give several applications of 
the theorems we have proved in sections 6 and 7.

The first special case of our theorems in Section 5, 6, and 7 is 
to $X=\RR^N$.

\begin{corollary}\sl Let $X=R^N$ with the standard
Euclidean metric and Lebesgue measure;
and let $K$ be a Calder\`on-Zygmund kernel such that $T_K$ is bounded on $L^2$.
 Then

(a) If $f\in BMO(\RR^N)$,
 then $C_f$ is bounded on $L^p(\RR^N)$ for all $1<p<\infty$;

(b) If $f\in BMO(\RR^N)$,
 then $C_f$ is bounded from $H^1(\RR^N)$ to $L^1_{\rm loc}(\RR^N)$;

(c) If $f\in VMO(\RR^N)$,
 then $C_f$ is compact on $L^p(\RR^N)$for all $1<p<\infty$;

(d) If $f\in \Lambda_{N(1/p-1)}(\RR^N)$ and $K$ is a $p$-standard
kernel, then $C_f: H^p(\RR^N)\to L^p(\RR^N)$
is bounded for all  $0<p<1$.
\end{corollary}

The result (a) generalizes the sufficient conditions in the main
theorem proved by Coifman, Rochbergh and Weiss in [CRW] and S. Janson
in [JAN]. The result (c) is a generalization of the main theorem
proved by Uchiyama in [UCH].

The converse of the above corollary is also true for 
some special Calder\`on-Zygmund kernels
(see [CRW], [JAN] and [UCH] for details.)

Next  we shall
 consider applications of the theorems we obtained in Sections 5, 6 and 7
to $X$ the boundary of a strictly pseudoconvex 
 domain  in $\CC^n$, a pseudoconvex
domain of finite type in $\CC^2$, or a convex domain of finite type in $\CC^n$.

Let $\Omega$ be a bounded domain in $\CC^n$ with $C^2$ boundary $\bomega$.
Let $d\sigma$ denote the Lebesgue surface measure on $\bomega$
and $L^p(\bomega)$
the usual Lebesgue space on $\bomega$ with respect to the measure $d\sigma$.
Let ${\H}^p(\Omega)$ be the holomorphic Hardy spaces defined in 
[KRA3] or [STE2].
Fatou's theorem [KRA1] shows that, for any $0<p\le \infty$, 
a holomorphic function
$f\in {\H}^p(\Omega)$ has a radial limit at almost all points on $\bomega$.
Thus one can identify ${\H}^p(\Omega)$ as a closed subspace 
of $L^p(\bomega)$.  Let $S: L^2(\bomega)\to {\H}^2(\Omega)$ be the
orthogonal projection via the reproducing kernel $S(z,w)$---the 
Szeg\"o kernel. 
For many special instances and classes of $\Omega$, 
we may identify the operator $S$ as a singular integral operator
on $\bomega$; in fact, in many instances
 $S(z,w)$ is given by a kernel that is $C^\infty$ 
on $\bomega\times \bomega\setminus D$ (where
$D$ is the diagonal of $\bomega\times \bomega$). 

First, we consider the case in which $\Omega$ is a strictly 
pseudoconvex domain
in $\CC^n$ with smooth boundary. Let $X=\bomega$. We denote by $d$ 
the usual quasi-metric on $\bomega$ defined in [STE2] or [FEF] or
[KRL1].  (In general, the quasi-metric defined on
 $\bomega$ is not symmetric, but we can define $\hat{d}(x,
y)=(1/2)(d(x,y)+d(y,x))$. Then $\hat{d}$ is a symmetric quasi-metric
having the same properties as $d$.) Also, we let the measure 
$d\mu=d\sigma$ be the usual Lebesgue/Hausdorff surface measure on $\bomega$. 
Then we have $|B(z,\delta)|=\sigma(B(z,\delta)\approx \delta^{n}$. It
is clear that condition (7.3) is satisfied with this measure. By
theorems in [BFG] and [BMS], the Szeg\"{o} kernel $S(z,w)\in
C^{\infty}(\bomega\times \bomega\setminus \{(x, x): x\in \bomega\})$
is a $p$-standard kernel (see the formulation in [KRL1]). 

The first main purpose of this section is to prove the following
theorem.

\begin{theorem}\sl  Let $\Omega$ be a bounded strictly pseudoconvex
domain in $\CC^n$ with smooth boundary.  With the notation above,
let $T_K=T_S$ ($S$ is the Szeg\"{o} kernel)
and let $f\in L^1(\bomeg)$. Then

(i) $f\in BMO(\bomeg)$ if and only if $C_f: L^p(\bomeg)\to
L^p(\bomeg)$ is bounded.

(ii) $f\in VMO(\bomeg)$ if and only if 
$C_f:L^p(\bomeg)\to L^p(\bomeg)$ is compact.
\end{theorem}

\proof Since $S(z,w)$ is a standard kernel on $\bomega\times \bomega$
and $T_S$ as a singular integral is the main part of the Szeg\"o projection (
In fact, one has $S(f)(z)=1/2 f(z)+c T_S(f)$, a.e. $z\in \bomega$.) it
is bounded on $L^2(\bomega)$. Therefore, the sufficient conditions 
for the commutator $C_f$ to be bounded and compact on
$L^p(\bomeg)$ for all $1<p<\infty$ in Theorem 8.2
follow directly from  Theorems 6.1 and 7.1.
 
Now we  prove the
necessary conditions of Theorem 8.2. 
\medskip

We first consider part (a).
 We assume that $f\in L^{p_0}(\bomeg)$ for some $1<p_0<<p$
and $C_f$ is bounded on $L^p(\Omega)$ for some $1<p<\infty$. We shall
show $f\in
BMO(\bomeg)$. Following [FEF],
[BFG] or [BMS], we let $\rho(z)$ be a strictly pluri-superharmonic 
defining function
for $\Omega$, and we set
$$
\psi(z, w)=\sum_{j=1}^n {\d \rho\over \d w_j}(w) (z_j-w_j)+(1/2)\sum_{jk}
{\d^2 \rho\over \d w_j\d w_k}(w) (z_j-w_j)(w_k-w_k) .
$$
Then there is a positive number $\delta>0$ such that
$$
S(z,w)=F(z,w) \psi(z,w)^{-n}+G(z,w) \log\psi(z,w)\leqno(8.3)
$$
for all $(z,w)\in R_\delta=\{(z,w)\in\bomega\times \bomega: d(z,w)<\delta\}$,
where $F, G\in C^\infty(\bomega\times \bomega)$ and $F(z,z)>0$ on $\bomega$. Moreover,
we have the following lemma proved by the authors [Lemma 5.2, KRL1].

\begin{lemma}\sl Let $\Omega$ be a smoothly bounded strictly pseudoconvex domain
in $\CC^n$. Then for each point $z_0\in \bomega$, there are holomorphic
functions $g_j$, and $C^\infty$ functions $h_j$ $(j=1,\cdots, M)$ and a function
$E(z,w)$ which is holomorphic in $z$ and $C^{n-1}(R_\delta)$ in $w$ such that
the reciprocal of the Szeg\"o kernel has the following decomposition property
(for $z$ close to $w$):
$$
{1\over S(z,w)}=\sum_{j=1}^M g_j(z, z_0) h_j(w,z_0) +E(z,w)\leqno(8.4)
$$
with
$$
|g_j(z, z_0)|\le C d(z, z_0)^{\gamma_j},\quad |h_j(z, z_0)|
\le Cd(z, z_0)^{\eta_j}\leqno(8.5)
$$
and $\gamma_j+\eta_j\ge n$ and $\eta_j,\ \gamma_j\ge 0$ integers; and
$$
|E(z,w)|\le C_M (|z-w|^M+|\psi(z,w)|^{2n-\epsilon})\leqno(8.6)
$$
for all $z, w\in B(z_0, \delta)$. 
\end{lemma}

Let $B=B(z_0,\delta)$ be any ball in $\bomeg$. If we 
choose $0<\epsilon\le n/p_0'$,
then for each $z\in B(z_0,\delta)$ we have
$$
\Big\{\int_B |\psi(z,w)|^{(n-\epsilon)p'}d\sigma(w)\Big\}^{1/p'}\le C\delta^{n}.
$$
Thus we have
\begin{eqnarray*}
\lefteqn{\int_B |f(z)-f_B| d\sigma(z)}\\
&=&{1\over |B|}\int_B \Big|\int_B (f(w)-f(z)) d\sigma(w)\Big| d\sigma(z)\\
&=&{1\over |B|}\int_B \Big|\int_B (f(w)-f(z)) S(z,w) S(z,w)^{-1} d\sigma(w)\Big|
d\sigma(z)\\
&=&{1\over |B|}\int_B \Big|\int_B (f(w)-f(z))\Big\{\sum_{j=1}^M g_j(z_0,z)
h_j(z_0,w)\\
&&\quad+ O(|z-w|^{2n}\delta^{n}+|\psi(z,w)|^{2n-\epsilon})
d\sigma(w)\Big\} \Big|d\sigma(z)\\
&\le& {1\over |B|}\int_B |\sum_{j=1}^M C_f(h_j(z_0,\cdot)
{\X}_B)(z) g_j(z_0,z) |d\sigma(z)\\
&&+{C\over |B|}\int_{B}|\int_{B}|f(w)-f(z)| |S(z,w)||z-w|^{2n}\delta^{n}
d\sigma(w) d\sigma(z)\\
&&+{C\over |B|}\int_B \int_B |f(w)-f(z)||\psi(z,w)|^{n-\epsilon}
d\sigma(w)d \sigma(z)\\
 &\le& {1\over |B|}\sum_{j+1}^M
||C_f(h_j(z_0,\cdot){\X}_B)||_p||g_j(z_0,\cdot){\X}_B||_{p'} 
+C ||f||_1 |B|\\
&\le& C{1\over |B|}\sum_{j=1}^MC||C_f|| ||h_j(z_0, \cdot) {\X}_B||_p ||g_j(z_0,\cdot)
{\X}_B)||_{p'}+C||f||_{p_0} |B|\\
&\le& {1\over |B|}CM||C_f|| \delta^{\beta_j} |B|^{1/p} \delta^{\alpha_j}
|B|^{1/p'}+C||f||_{p_0} |B|\\
&\le& CM||C_f||\delta^{n} +C||f||_{p_0}|B|\\
&\le& C(||f||_{p_0}+||C_f||)|B|.
\end{eqnarray*}
Therefore 
$$
{1\over |B|}\int_B |f(z)-f_B| d\sigma(z)\le C(||f||_{p_0}+||C_f||).
$$
for any ball $B$ in $\bomeg$.
This proves that $f\in BMO(\bomeg)$, and $||f||_{*}\le
C(||f||_{p_0}+||C_f||).$
\medskip

Next we prove part (b).
 Assume that $f\in L^{p_0}(\bomeg)$ with $1<p_0<<p<\infty$ and  
suppose that $C_f$ is
compact on $L^p(\bomeg)$ for some $p$.  We will 
show that $f\in VMO(\bomeg)$.

Since $C_f$ is compact, we have $C_f$ is bounded on $L^p(\bomeg)$. Thus we
have $f\in BMO(\bomeg)$. Next we show that
$$
{1\over |B(z,\delta)|}\int_{B(z,\delta)} |f(w)-f_B|d\sigma(w) \to 0
\leqno(8.7)
$$
uniformly for all $z\in \bomeg$ as $\delta\to 0$.

Suppose (8.7) is not true.  Then there are a sequence ${z_k}\subset \bomeg$
and ${\delta_k}$ and $\eta_0>0$ such that $\delta_k\to 0$ and 
$$
{1\over |B_k|}\int_{B_k} |f(w)-f_B|d\sigma(w)\ge \eta_0>0
\leqno(8.8)
$$
for all $k$ where $B_k=B(z_k,\delta_k)$.

For each $1\le j\le M$ and positive integer $k$, we let
 $$
\phi_{k,j}(w)={1\over |B_k|}h_j(z_k, w){\X}_{B_k} .
$$ 
Using the same argument as we did a few moments ago, we have
\begin{eqnarray*}
\lefteqn{\int_{B_k}|f(z)-f_{B_k}| d\sigma(z)}\\
&\le& C\int_{B_k} \Big| \sum_{j=1}^M C_f(\phi_{k,j})(z) g_j(z_k, z) \Big|d\sigma(z)\\
&&\quad + C\delta_k^{n-\epsilon} \int_{B_k}|f(z)-f_{B_k}|d\sigma(z).
\end{eqnarray*}
Therefore, when $k$ is big enough, we have
\begin{eqnarray*}
\lefteqn{\int_{B_k}|f(z)-f_{B_k}| d\sigma(z)}\\
&\le& C\int_{B_k} \Big| \sum_{j=1}^M C_f(\phi_{k,j})(z) g_j(z_k, z) \Big|d\sigma(z)\\
&=& C\int_{B_k} \Big| \sum_{j=1}^M C_f(\delta_k^{n/ p'-\eta_j}\phi_{k,j})(z) 
\delta_k^{-n/p'+\eta_j} g_j(z_k, z) \Big|d\sigma(z).
\end{eqnarray*}

Now we let
$$
y_{j,k}(z)=\delta_k^{n/p'-\eta_j}\phi_{k,j}(z)
$$
Then
$$
\|y_{j,k}\|_{L^p}\le C
$$
and
 $y_{k,j}\to 0$ weakly on $L^p(\bomeg)$ as $k\to \infty$ for all $1\le
j\le M$. In fact, for any smooth function $y(z)$, we have
\begin{eqnarray*}
\lefteqn{\int_{\bomega} y_{jk}(z) y(z) d\sigma(z)}\\
&=&{1\over |B_k|}\int_{B_k} \delta_k^{n/ p'-\eta_j} h_j(z_k, z) y(z) d\sigma(z)
\to 0 
\end{eqnarray*}
 as $ k\to \infty.$

Since $C_f$ is compact on $L^p(\bomeg)$, we have
$||C_f(y_{k,j})||_p\to 0$ as $k\to \infty$ for all $1\le j\le M$.
This implies that 
\begin{eqnarray*}
\lefteqn{{1\over |B(z_k,\delta_k)|}\int_{B(z_k,\delta_k)}
|f(w)-f_{B(z_k,\delta_k)}|d\sigma(w)}\\
&\le& C\|C_f(y_{j,k})\|_{L^p}
 \| \delta_k ^{ -n/p'+\eta_j} g_j{1\over |B_k|} {\X}_{B_k}\|_{p'}\\
&\le &C\|C_f(y_{j,k})\|_{L^p}
 \| \delta_k ^{ -n/p'+\gamma_j+\eta_j}{1\over |B_k|} \|B_k|^{1/p'}\\
&=&C\|C_f(y_{j,k})\|_{L^p}
 \| \delta_k ^{ -n/p'+\gamma_j+\eta_j-n+n/p'}\\
 &\le & C\|C_f(y_{j,k})\|_{L^p} \to 0.
\end{eqnarray*}
as $k\to \infty.$  This assertion contradicts (8.8).  
Therefore we have proved that $f\in
VMO(\bomeg)$. Part (b)  follows.

The proof of Theorem 8.2 is complete. \endpf

\begin{corollary}\sl Let $X=S^{2n-1}$, the unit sphere in 
$\CC^n$ or $\HH_n$, and
let $K=S$ be the corresponding Cauchy-Szeg\"o kernel. 
Let $f\in L^{p_0}(X)$ for some
$1<p_0<\infty$. Then

(i) $C_f: L^p(X)\to L^p(X)$ is bounded if and only if $f\in BMO(X)$ for $1<p<\infty$;

(ii) $C_f: L^p(X)\to L^p(X)$ is compact if and only if $f\in VMO(X)$ for $1<p<\infty$;

(iii) If $f\in BOM(S)$ then $C_f: H^p(S)\to L^1(S)$ is bounded;

(iv) If $0 < p <1$ and $f\in \Lambda_{n(1/p-1)}(S)$, 
the non-isotropic Zygmund space, then
$C_f: H^p(S)\to L^p(S)$ is bounded.
\end{corollary}

Part (i) was proved by Feldman and Rochberg in [FER] and part (iv) 
is related to some results in [LI].

Next we shall consider the case when $X=\bomega$,        
$\Omega$ a pseudoconvex domain
of finite type in $\CC^2$  or a convex domain of finite type in $\CC^n$.

When $\Omega$ is a smoothly bounded pseudoconvex domain of finite 
type in $\CC^2$,
then we shall use a variant of the quasi-metric defined in [NRSW]
as formulated in [KRL1] or [MCN1]. 
When $\Omega$ is a smoothly 
bounded convex domain in
$\CC^n$, we shall use the quasi-metric introduced by McNeal, 
for example see [MCN1],
[MCS1] and [KRL3].
 We shall
prove the following theorem:

\begin{theorem}\sl Let $\Omega$ be either a smoothly  bounded 
pseudoconvex domain of finite type in $\CC^2$ or a smoothly  bounded 
convex domain of finite type in $\CC^n$. Let $K=S$ be the Szeg\"o kernel
for $\Omega$. Let $f\in L^{p_0}(\bomega)$ for some $1<p_0<\infty$. Then
\smallskip \\
 
(i) $C_f: L^p(X)\to L^p(X)$ is bounded if $f\in BMO(X)$ for $1<p<\infty$;

(ii) $C_f: L^p(X)\to L^p(X)$ is compact  if $f\in VMO(X)$ for $1<p<\infty$.
\end{theorem}

In order to prove Theorem 8.5, we need some
more notation and some lemmas.
 Let $-r$ be a smooth defining function for $\Omega$ such that, in a
small neighborhood of $\bomeg$, $r(z)$ is the distance from $z$
to $\bomeg$. Let
$\nu(z)$ denote the unit inward normal vector to $\bomeg$ at $z\in \bomeg$. 
Then we may choose an $\epsilon_0>0$ small enough such that 
$$
z=\pi(z)+
r(z) \nu(\pi(z)),\quad \nu(z)= \left ({\d r\over \d \zbar_1}, \cdots, {\d
r\over \d \zbar_n}, {\d r\over \d z_1}, \cdots, {\d r\over \d z_n} \right )
$$
for all $z\in \Omega_{\epsilon_0}$.
 
 Let $L^2(\Omega)$ be the Lebesgue space on $\Omega$ and $A^2(\Omega)$ 
the subspace of holomorphic functions; 
we call $A^2$ the Bergman space. Let $P:L^2(\Omega)\to A^2(\Omega)$ 
be the orthogonal projection
---the Bergman projection with reproducing kernel $K(z,w)$, the 
Bergman kernel.

For each $ a\in
C^{\infty}(\Omegabar_{\epsilon_0})$, we define kernels $C(z,w)$ and
$C_{\epsilon}(z,w)$ on $\bomeg\times \bomeg$ by
$$
C(z,w)=\int_0^{\epsilon_0} a(z+t \nu(z)) K(z+t \nu(z), w) dt,
\leqno(8.9)
$$
and
$$
C_{\epsilon} (z,w)=\int_{\epsilon}^{\epsilon_0} a(z+t \nu(z)) K(z+t \nu(z), w) .
dt
\leqno(8.10)
$$

Following the main estimate in [NRSW] on $K(z,w)$ when $\Omega$ is a
smoothly bounded pseudoconvex domain of finite type in $\CC^2$ and
the estimate on $K(z,w)$ given in [MAN] and [MAS], and also
formulations in [KRL2] and [KRL3], we have

\begin{lemma}\sl  Let $\Omega$ be either a smoothly bounded
pseudoconvex domain of finite type in $\CC^2$ or
a smoothly bounded convex domain in $\CC^n$. Then
$C$ and $C_{\epsilon}$  as defined above are standard 
kernels on $\bomega\times \bomega$
\end{lemma}

Let $\Omega$ be  either a smoothly bounded
pseudoconvex domain of finite type in $\CC^2$ or
a smoothly bounded convex domain in $\CC^n$. Then, by the
results in [NRSW] for finite type domain in $\CC^2$, and those
in [MCN3] and [MCS1] for
convex domain of finite type in $\CC^n$, 
we have $P:L^p(\bomeg)\to {\H}^p(\Omega)$ is bounded for all
$1<p<\infty$. Let $\epsilon_0>$ be as the above. Then we
define an operator:
$$
Q(f)(z)=(P(f)(z+\epsilon_0 \nu(z)), \quad z\in \bomeg.
$$ 
It is easy to see that $Q:L^p(\Omega) \to C^{k}(\bomeg)$ is bounded for all
$1<p<\infty$ and $0\le k<\infty$. 

Then we have the following lemma proved in [KRL1] and [KRL3].

\begin{lemma}\sl Let $\Omega$ be either a smoothly bounded
pseudoconvex domain of finite type in $\CC^2$ or a
smoothly bounded convex domain in $\CC^n$. Then
$$ 
S(f)(z)=A(f)(z)+E S(f)(z)+Q(f),\quad z\in \Omega,\quad f\in L^p(\bomeg) ,
$$
where $A=\sum_{j=1}^n r_j I_{C_j}$, and $E=\sum_{j=1}^n [I_{C_j}, M_{r_j}]$
and
 $C_j(z,w) $ is obtained from (5.1) by replacing $a$ by $a_j$, $ j=1,\cdots,n.$
\end{lemma}

\begin{lemma} \sl Let $\Omega$ be either a smoothly bounded
pseudoconvex domain of finite type in $\CC^2$ or a
smoothly bounded convex domain in $\CC^n$.
 Let $f\in BMO(\bomeg)$ and $g\in C^1(\bomeg)$. Then $fg \in BMO(\bomeg)$.
\end{lemma}

\proof Let $B=B(z_0, \delta)$ be any ball on $\bomeg$. Then we have
$$
|g(z)-g(w)|\le ||g||_{C^1}|z-w|\le C \delta^{\eta}.
$$
for all $z, w \in B$ and for some $\eta>0$, where $\eta$ is a constant
depending only on the type of $\Omega$. Therefore
\begin{eqnarray*}
\lefteqn{|f(z)g(z)-(fg)_B|}\\
&=& |(f(z)-f_B)g(z)+ f_B g(z)-(fg)_B|\\
&=&|(f(z)-f_B)g(z) + {1\over |B|} \int_B (f(w) g(z)-f(w)g(w)) d\sigma(w)|\\
&\le& |f(z)-f_B||g(z)|+{1\over |B|} \int_B |f(w)||g(w)-g(z)| d\sigma(w)\\
&\le& |f(z)-f_B||g(z)|+C{1\over |B|} \int_B |f(w)| \delta^{\eta} d\sigma(w)\\
&\le & |f(z)-f_B||g(z)|+C \delta^{\eta} |f_B|\\
&\le& |f(z)-f_B||g(z)|+C \delta^{\eta} ||f||_* \log{1\over |B|}\\
&\le &|f(z)-f_B||g(z)|+C_{\eta} ||f||_*.
\end{eqnarray*}
Thus
$$
{1\over |B|} \int_B |fg(z)-(fg)_B|d\sigma(z)|\le C||g||_{C^1} ||f||_*.
$$
So we have
$$
||fg||_*\le C||f||_* ||g||_{C^1}
$$
and the proof of Lemma 8.8 is complete. \endpf

\noindent {\bf Remark:}  It should be noted that Stegenga [ST1] has given a
characterization, in the real variable setting, of those functions $g$
as in the last lemma that are $BMO$ multipliers.
\smallskip \\

 Now we apply Theorems 6.1, 7.1 and Lemmas 8.7 and 8.8 to see that
$$
[A, M_f]=\sum_{j=1} [r_j C_j, M_f]=[C_j, M_{fr_j}]
$$  
is bounded on
$ L^p(\bomeg)$ if $f\in BMO(\bomeg)$; and
$[A, M_f]$ is compact on $L^p(\bomega)$ if $f\in
VMO(\bomeg)$ for all
$1<p<\infty$. [We use here the fact that
$E$ and $Q$ are compact operators on $L^p$ for all
$1<p<\infty$.] 

Combining the above facts, the proof  (i) and (ii) of Theorem 8.5 is complete. 
\endpf 

Finally, we make the following remark. One may prove a result 
similar to (iii) in
Corollary 8.4 for a smoothly bounded pseudoconvex domain
of finite type in $\CC^2$ (with a different argument from the one in [KRL1])
or for a smoothly bounded convex domain
in $\CC^n$; one obtains these results as an 
application of Theorem 6.10.


\newpage

\begin{thebibliography}{ABCDE}

\bibitem[AIM]{abcd} H. Aimar and R. A. Mac\'ias, Weighted norm inequalities for the
Hardy-Littlewood maximal operators on space of Homogeneous type, {\it Proc. AMS}
 91(1984), 213--217.

\bibitem[BEL]{abcd} F. Beatrous and S-Y. Li, Boundedness and compactness of operators
of Hankel type, {\it  J. Functional Analysis}, 111(1993), 350--379.


\bibitem[BMS]{abcd} Boutet de Monvel and Sj\'ostrand, Sur la singularit\'e
des noyaux de Bergman et Szeg\"o, {\em Soc. Math. de France 
Asterisque}, 43--35(1976), 123--164.

\bibitem[CAL]{abcd} A. P. Calder\'on, Inequalities for maximal function relative to a
metric, {\it Studia Math.} 57 (1976), 297--306

\bibitem[CAZ]{abcd} A. P. Calder\'on and A. Zygmund, 
 On the existence of certain singular integrals, 
{\it Amer. J. Math.} 78(1956), 310--320.

\bibitem[CHR1]{abcd} M. Christ, A $T(b)$ theorem with remarks
on analytic capacity and the Cauchy integral, {\em Colloquium
Mathematicum} 50(1990), 601-628.

\bibitem[CHR2]{abcd} M. Christ, {\em Lectures on Singular Integral 
Operators},
 {\em CBMS Regional Conference Series in Mathematics},
number 77, 1990.

\bibitem[CHF]{abcd} M. Christ and R. Fefferman, A note on weighted
norm inequalities for the Hardy-Littlewood maximal operator,
{\em Proc. A.M.S.} 87(1983), 447-448.

\bibitem[COI]{abcd}  R. R. Coifman, A real variable characterization of $H^p$
{\em Studia Math.} 51(1974), 269--274.


\bibitem[COG]{abcd}  R. R. Coifman and L. Grafakos, 
Hardy space estimates for multilinear operators,
to appear. 

\bibitem[COW1]{abcd}  R. R. Coifman and G. Weiss, {\em Analyse Harmonique
Non-Commutative sur Certains Espaces Homogenes},
Springer Lecture Notes, vol. 242, Springer Verlag, Berlin, 1971.

\bibitem[COW2]{abcd}  R. R. Coifman and G. Weiss, Extensions
of Hardy spaces and
their use in analysis, {\em Bull. A. M. S.} 83(1977), 569-643.

\bibitem[CRW]{abcd} 
  R. R. Coifman, R. Rochberg and G. Weiss, Factorization theorems for
Hardy spaces in several variables, {\it Ann. Math.} 103(1976), 611-635.

\bibitem[CIS]{abcd} 
  J. Cima and D. Stegenga, Hankel operators on $H^p$
domain in $\CC^n$, {\it Analysis at Urbana 1  LMS}, Lecture Note 137, 133-150.

\bibitem[FEF]{abcd} 
 C. Fefferman, The Bergman kernel and biholomorphic mapping of pseudoconvex domains,
{\em Invention Math.} 26(1974), 1--65.

\bibitem[FER]{abcd} 
 M. Feldman and R. Rochberg, Singular value estimates for commutators
and Hankel operators on the unit ball and the Heisenberg group, {\it Analysis and
P.D.E.} Lecture Notes in P.A. Math. 122, Dekker, New York, 1990.

\bibitem[FES]{abcd}  C. Fefferman and E. M. Stein, $H^p$ spaces of several
variables, {\em Acta Math.} 129(1972), 137-193.

\bibitem[FEL]{abcd}  M. Feldman, Mean Oscillation, Weighted Bergman spaces
and Besov spaces on the Heisenberg group and
atomic decomposition, preprint.

\bibitem[GAL]{abcd} J. Garnett and  R. Latter, 
The atomic decomposition for Hardy spaces
in several complex variables, {\em Duke Math. J. } 45 (1978), 815--845.

\bibitem[GL]{abcd} J. Grafakos and  X. Li, preprint.


\bibitem[HAW]{abcd} Y-S. Han and  G. Weiss, 
Function spaces on
spaces of homogeneous type, preprint.

\bibitem[JAN]{abcd}  S. Janson, Mean oscillation and commutators of singular integral
operators, {\it Ark. Math. J.} (1978).

\bibitem[KRA1]{abcd}  S. G. Krantz, Geometric Lipschitz 
spaces and applications
to complex function theory and nilpotent groups,
 {\em J. Functional Analysis} 34(1979), 456--471.


\bibitem[KRA2]{abcd}  S. G. Krantz, Geometric Analysis and
  Function Spaces, {\em CBMS} 
Regional Conference Series in Mathematics, No 81 (1993).


\bibitem[KRL1]{abcd} S. G. Krantz and  S-Y.  Li, 
On the decomposition theorems for Hardy spaces and applications
in domains in $\CC^n$, {\em J. of Fourier Analysis}, to appear.

\bibitem[KRL2]{abcd} S. G. Krantz and  S-Y.  Li, 
A Note  on Hardy spaces and functions of bounded mean oscillation on 
domains in $\CC^n$, {\em Michigan J. of Math}, 1994.

\bibitem[KRL3]{abcd} S. G. Krantz and  S-Y.  Li, 
Duality  theorems for Hardy and Bergman spaces on convex domains of
finite type in $\CC^n$, {\em Ann. of Fourier Institute}, to appear.


\bibitem[LAT]{abcd} R. H. Latter, A characterization of $H^p(\RR^N)$ in terms of atoms,
{\em Studia Math.} 62(1978), 92--101.

\bibitem[LI]{abcd} S-Y. Li, Toeplitz Operators on Hardy Space $H^p(S)$ with $0<p\le 1$,
{\em Integral Equations and Operator Theory}, 15(1992), 808--824.

\bibitem[MS1]{abcd}  R. Macias and C. Segovia, Lipschitz functions on spaces
of homogeneous type, {\em Adv. in Math.} 33(1979), 257-270.

\bibitem[MS2]{abcd}  R. Macias and C. Segovia,
A decomposition into atoms of distributions on spaces
of homogeneous type, {\em Adv. in Math.} 33(1979), 271-309.


\bibitem[MCN]{abcd} J. McNeal, Estimates on the Bergman kernels of convex
 domains,
{\it Advances in Math.}, in press.

\bibitem[MCS]{abcd} J. D. McNeal and E. M. Stein, Mapping properties
of the Bergman projection on convex domains of finite type, preprint.


\bibitem[NRSW]{abcd} A. Nagel, Rosay, E. M. Stein and Wainger,  Estimates
for the Bergman and Szeg\"o kernels in $\CC^2$,
{\em Ann. of Math.},
129(1989), 113--149.

\bibitem[ROS]{abcd}  L. R. Rothschild and E. M. Stein,
Hypoelliptic differential operators and nilpotent groups,
 {\em Acta Math.} 131(1976), 247-320.


\bibitem[ST1]{abcd} D. A. Stegenga, Bounded
 Toeplitz operators on $H^1$ and applications of duality between $H^1$ and 
functions of bounded mean oscillation, {\em Amer J. of Math.} (1976), 573--589.

\bibitem[STE1]{abcd}  E. M. Stein,
Harmonic Analysis,
 {\em Princeton University Press}, Princeton, New Jersey, 1993.

\bibitem[STE2]{abcd} E. M. Stein, {\em Boundary
Behavior of Holomorphic Functions of Several
Complex Variables}, Princeton University Press,
Princeton, 1972.

\bibitem[STT]{abcd}J. O. Str\"omberg and A. Torchinski, Weights, sharp maximal
functions and Hardy spaces, {\it Bull. Amer. Math. Sco.} 3(1980), 1053--1056.

\bibitem[UCH]{abcd} A. Uchiyama, Compactness of operators of Hankel type, {\em T\^ohoku
Math. J. } 30(1978), 163--171.
\end{thebibliography}
\end{document}